\documentclass{article}
\usepackage{amscd,amsmath,amssymb,amsfonts,emlines}
\usepackage[dvips]{graphicx}
\makeatletter \@addtoreset{equation}{section} \makeatother
\def\Z{\mathbb{Z}} \def\R{\mathbb{R}}
\def\C{\mathbb{C}} \def\H{\mathcal{H}}
\def\N{\mathbb{N}}
\title{$h$-holomorphic functions of double variable\\and their applications}
\author{Dmitriy G. Pavlov\thanks{geom2004@mail.ru}\ \  {\em and}\ \  Sergey S. Kokarev\thanks{logos-center@mail.ru}}
\date{Research Institute of Hypercomplex system in Geometry and Physics (Fryazino), Regional Scientific-Educational Center
\, "Logos"\,, (Yaroslavl) }

\begin{document}\maketitle
\begin{abstract}
The paper studies the complex differentiable functions of double argument and their properties,
    which are similar to the properties of the holomorphic functions of complex variable: the
    Cauchy formula, the hyperbolic harmonicity, the properties of general $h$-conformal mappings
    and the properties of the mappings, which are hyperbolic analogues of complex elementary functions.
    We discuss the utility of $h$-conformal mappings to solving 2-dimensional hyperbolic
    problems of Mathematical Physics.
\end{abstract}
{\bf Keywords:} double numbers, hyperbolic Cauchy problem, $h$-holomorphicity, conformal mappings.
%
%
\section{Introduction}

Double numbers (or, as they are called sometimes, hyperbolic-complex or splittable numbers) are known for
    quite a long time, and find applications both in Mathematics and Physics \cite{1,2}. Taking into account
    that their corresponding algebra is isomorphic to the direct sum of two real algebras, it was considered
    for a long time that the properties of double numbers are not interesting, and moreover, that these are
    even trivial, while compared to complex numbers.\par
In this paper we try to prove the wrongness of this opinion and to show that the possibilities offered by
    the algebra and the analysis (i.e., the $h$-analysis) of double numbers, considered together with the
    associated geometry of the 2-dimensional plane, are still far from being exhaustively studied.\par
We shall see that, in many respects, double numbers are far -- in any sense, from being inferior to
    complex numbers. Considering the important circumstance that the geometry of double numbers is
    pseudo-Euclidean and hence corresponding to the 2-dimensional geometry of Space-Time, we obtain
    a new framework fruitful in ideas for Physics,  which have a purely algebraic essence.
A special particular feature of our study expresses the fact that the $h$-analytic mappings of double
    variable cover an infinite-dimensional space -- which is exactly the case for the functions of
    one complex variable.
Moreover, to each mapping defined on the complex plane, it corresponds a unique $h$-analytic mapping
    of double variable, and converse.\par
%
In Section \ref{doubleee}, we concisely present properties of the plane of double
    variable and define the hyperbolic polar coordinate system, cones, $h$-conformal mappings and study
    their general properties. In Section \ref{holooo}, we define the $h$-holomorphic mappings and discuss
    several of their properties, which include analogs of the Cauchy formula and the Cauchy theorem, of the
    Cauchy-Riemann conditions, and of $h$-conformal mappings. In Section \ref{elemmm}, we define and study
    the properties of standard hyperbolic elementary functions of double variable.
At last, in the Conclusions section, we discuss the potential applications and the perspectives of further
    work in developing the present theory of functions of double variable.
\section{Double numbers}\label{doubleee}

By analogy with the algebra of complex numbers $\C$, we define the algebra of double numbers $\H $ by means of
    two generators $\{1,j\}$ of the 2-dimensional $\R$-modules with the related multiplication table:
\begin{equation}\label{alg2}\begin{array}{c|c|c} & 1& j\\ \hline 1& 1 & j\\ \hline j& j& 1\end{array}.\end{equation}
Having in view the further application of this algebra to describing the 2-dimensional
    Space-Time, we shall denote the elements of $\H $ as: $h=1\cdot t+jx\in\H$, where $t,x\in\R$.
    Similarly to the complex numbers, the real number $\text{Re}\,h\equiv t$ will be called
    {\em the real part} of the double number $h$, and the real number $\text{Im}\,h\equiv x$
    will be called {\em the imaginary part } of the double number $h$. The algebra of double numbers
    with the multiplication table (\ref{alg2}) does not determine a numerical field, since it contains
    zero-divisors, i.e., the equation $h_1h_2=0$ may be satisfied by nonzero elements $h_1$ and $h_2$.
This is one of the reasons for which double numbers were not widely used in applications, as the
    complex numbers did. But exactly this feature shows in algebraic terms a very important
    circumstance of the 2-dimensional Space-Time -- the occurrence of light-cones.
    The geometrical interpretation of double numbers is analogous to the interpretation
    of complex numbers: in the plane of double variable (in brief -- on the hyperbolic plane),
    to each double number there corresponds a position-vector, whose coordinates are the
    real and the imaginary parts of this number. Here, the sum and the difference of double numbers
    is represented by the standard parallelogram rule for the corresponding position vectors
    on the hyperbolic plane

The involutive operation of {\em complex conjugation} for double numbers can be defined in
    the following way: $h=t+jx\mapsto \bar h=t-jx$.
Geometrically, this operation describes
    the reflection of the hyperbolic plane with respect to the axis $\text{Im}\, h=0$.
    Like in the complex case, the couple $(h,\bar h)$ can be regarded as independent double
    coordinates on the hyperbolic plane, which relate to the Cartesian coordinates by means
    of the formulas:
\begin{equation}\label{coordh}
x=\frac{h+\bar h}{2};\quad  y=\frac{h-\bar h}{2j}.\end{equation}

The complex coordinate bilinear form $\mathcal{G}\equiv dh\otimes
d\bar h$ splits   into
    its symmetric $\Xi$ and skew-symmetric $\Omega$ irreducible components, as follows:
\begin{equation}\label{2formh}\mathcal{B}\equiv dh\otimes d\bar h=\Xi-j\Omega,\end{equation}
where $\Xi=dt\otimes dt-dx\otimes dx$ is a pseudo-Euclidean metric form, and
    $\Omega\equiv dt\wedge dx=-j dh\wedge d\bar h/2$ is the 2-dimensional volume form.
We note, that the algebra of double numbers induces on the plane of double variable
    a 2-dimensional pseudo-Euclidean (hyperbolic) geometry endowed with the metric form
    $\Xi$, which justifies the established by us denomination "hyperbolic plane".

\begin{figure}[htb]\centering \unitlength=0.50mm \special{em:linewidth 0.4pt}
    \linethickness{0.4pt} \footnotesize \unitlength=0.70mm
    \special{em:linewidth 0.4pt} \linethickness{0.4pt}
    \input{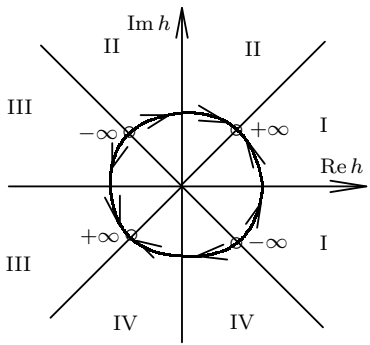}
    \caption{\small The domain $\R\sqcup-\R\sqcup\R\sqcup-\R$ of the change of the angle
    $\psi$ on the plane $\H$. The orientation is synchronized in opposed quadrants
    and is opposed in the neighboring ones. For different angles in different quadrants
    one can enumerate the angle $\psi$ by using the index $k$: $\psi_k$, $(k=1,2,3,4)$. }
    \label{psi}
\end{figure}

The transition to hyperbolic polar coordinates and to the exponential form of representing
    double numbers exhibit a series of special features, which are absent in the case of
    complex numbers. The pair of lines $t\pm x=0$ contains the subset of double numbers
    with zero squared norm\footnote{Strictly speaking, the occurrence of zero-divisors and
    the possibility of having negative values for the expression $h\bar h$ does not
    allow us to speak about the norm of a double number in its rigorous meaning}.
For the brevity of terms and for preserving the partial analogy to the complex numbers,
    we shall call the quantity $\sqrt{|h\bar h|}$ {\em the norm} or {\em the module}
    of the double number (see further formula \ref{regih}). The lines of double numbers whose
    square of the norm vanishes, split the whole hyperbolic plane into four quadrant-type domains,
    which are represented on the drawing by the numbers I, II, III and IV (Fig.\ref{psi}).

One can immediately verify that in each of the mentioned areas, the double numbers allow
    a hyperbolic polar representation of the form:
\begin{equation}\label{trigh}h=t+jx=\epsilon\varrho(\cosh\psi+j\sinh\psi),\end{equation}
where for each quadrant take place the following definitions of quantities:
\begin{equation}\label{regih}\begin{array}{rlcc}
    \text{I}:&\epsilon=1,& \varrho=\sqrt{t^2-x^2},\quad \psi=\text{Arth}(x/t);\\
    \text{II}:&\epsilon=j,& \varrho=\sqrt{x^2-t^2},\quad \psi=\text{Arth}(t/x);\\
    \text{III}:&\epsilon=-1,& \varrho=\sqrt{t^2-x^2},\quad \psi=\text{Arth}(x/t);\\
    \text{IV}:&\epsilon=-j,& \varrho=\sqrt{x^2-t^2},\quad \psi=\text{Arth}(t/x).
    \end{array}
\end{equation}
The quantities $\varrho$ and $\psi$, defined in each of the quadrants by the formulas
    (\ref{regih}), will be called {\em the module} and respectively {\em the argument}
    of the double number $h$.
In this way, in each of the quadrants we have $0\le\varrho<\infty$, and the
    quadrants themselves are parametrized by different copies of real lines,
    which together determine {\em the manifold $\Psi$ of angular variables}
    as a oriented disjoint sum $\R\sqcup -\R\sqcup \R\sqcup -\R$.
    Moreover, the manifold $\Psi$ can be suggestively represented by compactifying
    each copy of $\R$ into an open interval and further by gluing the intervals at their ends
    to obtain a circle with four pinched points.

\begin{figure}[htb]\centering\input{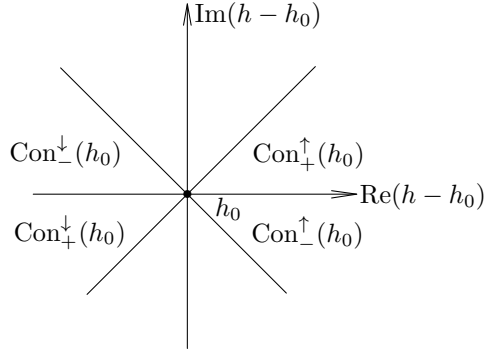}
    \caption{\small The definition of the cones and subcones of the point $h_0$.
    $\text{Con}(h_0)=\text{Con}_+(h_0)\cup\text{Con}_-(h_0)=(\text{Con}_+^\uparrow(h_0)
    \cup\text{Con}_+^\downarrow(h_0))\cup(\text{Con}_-^\uparrow(h_0)\cup\text{Con}_-^\downarrow(h_0))$.}
    \label{conus}
\end{figure}

We note, that the set of double numbers of zero norm is not described in any of the
    coordinate charts introduced above of the hyperbolic polar coordinate system.
    In the following, we shall call the subset of double numbers of the form
\begin{equation}\label{con}h_0+h(1\pm j),\end{equation}
(where $h_0$ --- some fixed double number, $h$ is any double
number) as {\em the cone of the number $h_0$}
    and will be denoted as $\text{Con}(h_0)$. All the points which lie in $\text{Con}(h_0)$,
    have their hyperbolic distance to the point $h_0$ equal to zero.
    Sometimes we shall make a difference between $\text{Con}_+(h_0)$ and
    $\text{Con}_-(h_0)$, depending accordingly on the signs in (\ref{con}).
    In the same way, we can distinguish the sub-cones $\text{Con}_+^\uparrow(h_0)$ and
    $\text{Con}_+^\downarrow(h_0)$, relative to the cases $\text{Re}h+\text{Im} h>0$
    and $\text{Re}h+\text{Im} h<0$ accordingly and the sub-cones
    $\text{Con}_-^\uparrow(h_0)$ and $\text{Con}_-^\downarrow(h_0)$ for the cases
    $\text{Re}h-\text{Im} h>0$ and $\text{Re}h-\text{Im} h<0$, accordingly.
    All the cones and subcones are shown in Fig.\ref{conus}.

The hyperbolic Euler formula: $\cosh\psi+j\sinh\psi=e^{j\psi}$ can be verified by
    expanding the left and right sides into formal Maclaurin series, and comparing
    their real and imaginary parts. The hyperbolic Euler formula leads to the exponential
    representation of double numbers:
\begin{equation}\label{exph}h=t+jx=\epsilon\varrho e^{j\psi}=\epsilon e^{\Theta},\end{equation}
where in the last equality we passed to the "complex hyperbolic angle"
\begin{equation}\label{lnh}\Theta=\ln\varrho+j\psi\equiv\ln
h,\end{equation}
having sense in the first quadrant.  So, the product of a pair of
double numbers reduces to adding their hyperbolic angles
    and the product of the sign factors $\epsilon$.

The formulas for computing lengths of curves and areas of domains in $\H$ reads as follows:
\begin{equation}\label{measc}
\text{Length}[\gamma]=\int\limits_{\tau_A}^{\tau_B}\sqrt{|\dot h\,
    \dot{\bar h}|}\, d\tau;\quad
    \text{Area}[\Sigma]=-\frac{j}{2}\int\limits_{\Sigma}dh
    \wedge d\bar h=-\frac{j}{4}\oint\limits_{\partial\Sigma}(h\, d\bar h-\bar h\, dh),\end{equation}
where in the last equality we used the double number variant of the foundamental Poincar\'e-Darboux
    Theorem regarding the integration of differential forms.

\medskip

{\small {\bf Example}. We shall compute the length of the arc of Euclidean circle of Euclidean radius $r$
    whose center is zero, encompassed between the points $1$ and $j$ on the plane of double variable.
Substituting in the Euclidean equation of the circle $t^2+x^2=r^2$ the polar hyperbolic coordinates:
    $t=\varrho\cosh\psi$, $x=\varrho\sinh\psi$, we get the hyperbolic polar equation of the Euclidean circle:
\begin{equation}\label{pcirc}\varrho(\psi)=\frac{r}{(\cosh^2\psi+\sinh^2\psi)^{1/2}}.\end{equation}
Constructing the pseudo-Euclidean line-element: $dl^2=|d\rho^2-\rho^2\,d\psi^2|$ along the circle,
    by considering (\ref{pcirc}), we get after several elementary calculations and differentiation:
\begin{equation}\label{dlina}dl=\frac{r d\psi}{(\cosh^2\psi+\sinh^2\psi)^{3/2}}\end{equation}
Due to the symmetry of the arc relative to the first bisector $t=x$, it suffices to compute the length of the
    half of the arc in which $\psi$ varies from $0$ to $\infty$, and then we double the obtained result.
The integral which provides this length is:
\[L=2\int\limits_0^{\infty}\frac{r d\psi}{(\cosh^2\psi+\sinh^2\psi)^{3/2}}.\]
Using the substitution:  $\tanh\psi=\xi$, this integral reduces to a simpler form, and can
    be expressed in terms of elliptic integrals of first and second order:
\[L=2r\int\limits_0^1\sqrt{\frac{1-\xi^2}{1+\xi^2}}\frac{d\xi}{1+\xi^2}=2\sqrt{2}
    r[E(1/\sqrt2)-K(1/\sqrt2)/2]\approx1.2r\]}
%
%
\section{$h$-holomorphic functions of double variable}\label{holooo}

The function $\ln h$, defined by formula (\ref{lnh}), is a simple
and important sample of the class of so-called {\em
$h$-holomorphic mappings of double variable}, whose definition
emerges from considerations similar to those which yield the
definition of holomorphic functions of complex variable. Any
smooth mapping $f:\R^2\to\R^2$ can be represented by a pair of
real components, and we can pass (with using of (\ref{coordh})) to its
representation by means of a pair of double variables $\{h,\bar
h\}$, as follows:
\begin{equation}\label{mapch}(h,\bar h)\mapsto(h',\bar h'):\ \
    h'=F_1(h,\bar h);\quad \bar h'=F_2(h,\bar h).\end{equation}
If $\R^2$ is regarded now as the plane $\H$ of double variable,
then
    we can naturally limit ourselves to the mappings which preserve
    the hyperbolic complex structure of the plane, i.e., such
    mappings of the form $h\in\H \to s=F(h)\in \H$.
The differentiable functions $\R^2\to\R^2$, which satisfy
    the condition%
\footnote{The notion of derivative of a function
    $F(h,\bar h)$ relative to its arguments is similar to the
    one used in real analysis. Namely, we can define the
    differentiability of a function $F$ at the point $(h,\bar h)$
    as the following property of its variation:
    $\Delta F=A(h,\bar h)\,\Delta h+B(h,\bar h)\,
    \Delta\bar h+o(\|\Delta h\|_H)$, where $\|\Delta h\|_H
    \equiv[\Delta t^2-\Delta x^2]^{1/2}$ is the pseudo-Euclidean
    norm of the variation of the variable.
Passing to different limits for $\|\Delta h\|_H\to 0$, we obtain the definition
    of partial derivatives or "directional derivatives". We encounter
    several significant problems, related to the definition of the
    convergence and the limit by natural for the double numbers hyperbolic
    norm. In this paper we shall not
    deal with these purely mathematical questions and we use only
    those operations and properties, whose definitions are clear, though
    to a certain extent, formal.}:
\begin{equation}\label{anh}F_{,\bar h}=0\end{equation}
are called {\em $h$-holomorphic mappings of double variable $h$}.
    The functions which satisfy the condition:
\begin{equation}\label{aanh}F_{,h}=0\end{equation}
are called {\em anti-holomorphic} mappings of double variable.

By analogy to holomorphic functions of complex variable, the
    holomorphic functions of double variable can be defined by formal
    power series, whose convergence often follows from the convergence
    of the corresponding real series.

\bigskip

{\small {\bf Example}. The following identities can
    be straightforward verified by means of expansion into formal
    series:
\begin{equation}\label{eqqh}S(jx)=jS(x);\quad C(jx)=C(x);\quad \end{equation}
\[S(h)=S(t+jx)=S(t)C(x)+jC(t)S(x);\quad C(h)=C(t+jx)=C(t)C(x)-jS(t)S(x),\]
where $x\in\R$, $S$ is the sinus (elliptic or hyperbolic), $C$ is
    the cosine (elliptic or hyperbolic) in the left and the right hand
    sides of the equalities, accordingly, which are defined by their
    standard series.

In fact, these equalities are particular cases of the more general
    identity:
\[f(jx)=S_f(x)+jA_f(x),\]
where $S_f\equiv [f(x)+f(-x)]/2$, $A_f\equiv [f(x)-f(-x)]/2$ are,
    respectively, the symmetric and the skew-symmetric parts of the
    arbitrary analytic function $f$.

We shall show, that a holomorphic function always maps zero-divisors
    to zero-divisors. The proof relies on the following formal identity:
\begin{equation}\label{eqqh1}(1\pm j)^{\alpha}\equiv 2^{\alpha-1}(1\pm j),\ \alpha\in R.\end{equation}
For $\alpha\in\N$, the identity immediately emerges from the simpler
    one: $(1\pm j)^2=2(1\pm j)$. For arbitrary $\alpha$, we need to use the
    expansion formulas into Maclaurin series. From one side, one has
\begin{equation}\label{row1}(1\pm j)^\alpha=1\pm\alpha j+\frac{\alpha(\alpha-1)}{2!}\pm
    \frac{\alpha(\alpha-1)(\alpha-2)}{3!}j+\frac{\alpha(\alpha-1)(\alpha-2)(\alpha-3)}{4!}+\dots
\end{equation}
On the other side,
\begin{equation}\label{row2}2^{\alpha-1}=(1+1)^{\alpha-1}=1+\alpha-1
    +\frac{(\alpha-1)(\alpha-2)}{2!}+\frac{(\alpha-1)(\alpha-2)(\alpha-3)}{3!}+
    \frac{(\alpha-1)(\alpha-2)(\alpha-3)(\alpha-4)}{4!}+\dots
\end{equation}
Multiplying this line, term by term, with $(1\pm j)$, we get:
\begin{equation}\label{row3}
    2^{\alpha-1}(1\pm j)=1\pm j+(\alpha-1)\pm(\alpha-1)j+
    \frac{(\alpha-1)(\alpha-2)}{2!}\pm\frac{(\alpha-1)(\alpha-2)}{2!}j+
\end{equation}
\[\frac{(\alpha-1)(\alpha-2)(\alpha-3)}{3!}\pm\frac{(\alpha-1)(\alpha-2)(\alpha-3)}{3!}j+\dots\]
By combining in the expansion (\ref{row3}) all the successive pairs
    of even terms (which contain $j$), we obtain all the successive even
    terms of the series (\ref{row1}), and combining in (\ref{row3}) all
    the successive pairs of even terms (beginning with the pair
    "third-fifth"), we get all the successive odd terms of the series
    (\ref{row1}) (beginning with the third one). In this way, the formal
    series from the left and right sides of (\ref{eqqh1}) coincide,
    q.e.d.

Considering now a holomorphic function $F(h)$ given by a power
    series:
\begin{equation}\label{row4}F(h)=\sum\limits_{k=0}^\infty c_k(h-h_0)^k,\end{equation}
due to the identity (\ref{eqqh1}), we get on the cone
    $\text{Con}(h_0)$ of the arbitrary point $h_0$:
\[F(h)|_{h\in\text{Con}(h_0)}=F(h_0)+\sum\limits_{k=1}^\infty c_k((t-t_0)\pm
    j(t-t_0))^k=F(h_0)+(1\pm j)\sum\limits_{k=0}^\infty
    c_k2^{k-1}(t-t_0)^k\subset\text{Con}(F(h_0)),
\]
which proves the claim. In fact, as we notice from the obtained
    expression, the holomorphic mapping can provide an inversion of the
    cone (i.e., it can maps the component $\text{Con}_{\pm}^\uparrow$
    into the corresponding component $\text{Con}_{\pm}^\downarrow$ and
    converse), but it cannot map its branches $\text{Con}_+$ and
    $\text{Con}_-$ one into another. It is easy to check that the last
    property is achieved by means of anti-holomorphic mappings}.
%
%
\subsection{Hyperbolic Cauchy-Riemann conditions}

Let's  write the condition (\ref{anh}) in Cartesian coordinates:
\[F_{,\bar h}=(U+jV)_{,\bar h}=\frac{(U+jV)_{,t}}{\bar
    h_{,t}}+\frac{(U+jV)_{,x}}{\bar h_{,x}}=U_{,t}-V_{,x}+j(V_{,t}-U_{,x})=0.
\]
This implies {\em the Cauchy-Riemann condition of hyperbolic
    analyticity}:
\begin{equation}\label{crh}U_{,t}=V_{,x};\quad U_{,x}=V_{,t}.\end{equation}
If the mapping $F$ is $h$-holomorphic in the sense of the former
    definition, then the mapping $\ln F=\ln\varrho_F+j\psi_F+
    \ln\epsilon_F$ is $h$-holomorphic as well. This leads to
    the Cauchy-Riemann condition, written in terms of module and
    argument of the function of double variable:%
\footnote{We need to say that, in general, $\epsilon_F$ is equal to
    a certain constant inside each quadrant on the image plane,
    which abruptly changes while passing over the borders of the
    quadrants. For this reason, the formulas (\ref{crh1}) are correctly
    defined on $\H $ minus the cross-shape lines of its
    zero-divisors. These considerations remove the question
    regarding the meaning of the expression $\ln\epsilon_F$, which may
    arise while going beyond the framework of the double numbers
    algebra.}
\begin{equation}\label{crh1}(\ln\varrho_F)_{,t}=(\psi_F)_{,x};\quad
    (\ln\varrho_F)_{,x}=(\psi_F)_{,t}.\end{equation}

It is easy to check that from the conditions (\ref{crh}), it follows
    the hyperbolic harmonicity of the real and of the imaginary parts of
    the holomorphic function $F$, which is expressed by the equations:
\begin{equation}\label{harmh}\Box U=\Box V=0,\end{equation}
where $\Box\equiv \partial^2_{t}-\partial^2_x$ is the wave operator of second order (the "hyperbolic Laplacian").
%
%
\subsection{The hyperbolic analogues of the Cauchy theorem}

 In order to stress the independence of our  proof of Cauchy theorem
    on the metrical properties of the double (or complex) plane,
    we shall provide it in terms of differential forms%
\footnote{We do not limit ourselves to defining curvilinear
    integrals on the plane of double variable, having in view their
    convergence to a pair of integrals of a 1-form on the Cartesian
    plane, whose definition is standard}.
For any simple closed contour $\Gamma$ which surrounds the domain
    $\Sigma\subset\H $ and for any holomorphic function of double variable
    $F=U+jV$, we have the following chain of equalities:
\[\oint\limits_{\Gamma}F(h)\, dh=\oint\limits_{\Gamma}U\, dt+V\,
    dx+j\oint\limits_{\Gamma}U\, dx+V\, dt=\int\limits_{\Sigma}
    [(V_{,t}-U_{,x})+j(U_{,t}-V_{,x})\,dt\wedge dx=0
\]
as consequence of the conditions (\ref{crh}). The second sign of the
    equality expresses the Poincar\'e-Darboux Theorem regarding
    integration of 1-forms along closed paths. Using the complex
    analysis terminology, the proof looks even shorter:
\[\oint\limits_{\Gamma}F(h)\, dh=\int\limits_{\Sigma}F_{\bar h}\,
    d\bar h\wedge dh=0\]
if taking into account (\ref{anh}). From purely topological
    considerations, similar to the ones from the complex plane,
    the integral of a holomorphic function cancels on the border of
    a multi-connected domain.

\begin{figure}[htb]\centering\input{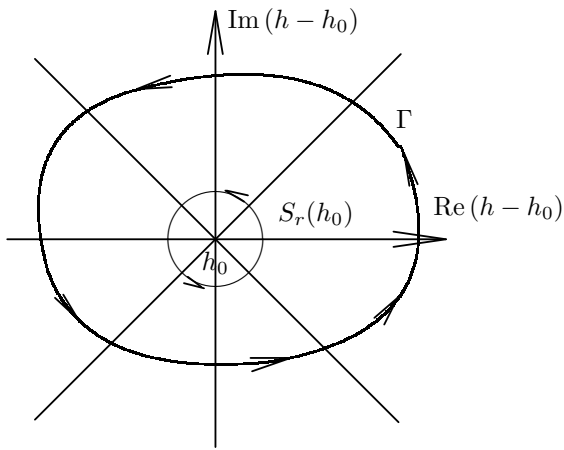}
    \caption{\small Towards setting the Cauchy Theorem on the plane of double variable}.\label{psi1}
\end{figure}

For the Cauchy integral in its hyperbolic version, we have now the
    equality:%
\footnote{We need to remark a certain ambiguity in the notation
    (\ref{intch}): the division in the integrand is meaningless
    on the intersection $\Gamma\cap\text{Con}(h_0)$. In this way,
    rigorously speaking, we need to remove the points $\text{Con}(h_0)$
    from the domain of the integrant, and the integral should be
    regarded as the limit of the integral over a non-connected
    contour, whose discontinuities are concentrated in the neighborhood
    $\text{Con}(h_0)$ and their (Euclidean) measure tends to zero.
Our results correspond to such an integral, considered in the sense
    of its principal value. For its existence, the contour has to
    transversally approach the lines of the cone. We do not insist on
    these purely mathematical questions in this paper, and postpone
    their more detailed investigation for a complementary paper.}
\begin{equation}\label{intch}\oint\limits_{\Gamma}\frac{F(h)}{h-h_0}\,
    dh=\oint\limits_{S_r(h_0)}\frac{F(h)}{h-h_0}\, dh,\end{equation}
which follows from the hyperbolic Cauchy Theorem. Here $S_r(h_0)$ is
    the (Euclidean) circle of radius $r$ and center at the point $h_0$,
    and the integral does not depend on the radius of this circle
    (see Fig.\ref{psi1}). We perform the change of variable:
    $h=h_0+\epsilon\varrho(r,\psi)e^{j\psi}$, where the mapping
    $\varrho(r,\psi)e^{j\psi}=rf(\psi)e^{j\psi}$ is the polar
    parametrization of the Euclidean circle $S_r(0)$ in terms of the
    hyperbolic polar coordinate system. We have, as follows from the
    example in Section 3, $f=1/\sqrt{\cosh^2\psi+\sinh^2\psi}$.
    We further need only the bijectivity of the function $f$.
    Integrating by $\psi$, we infer: $h-h_0=\epsilon rf(\psi)e^{j\psi}$, $
    dh=\epsilon r(df+jf\, d\psi)e^{j\psi}$, and the Cauchy integral
    gets the form:
\[\oint\limits_{S_r(h_0)}F(h)(d\ln f+j\,d\psi).\]
Using the $r$-independence of the integral and passing to limit for
    $r\to 0$, we obtain:
\[\oint\limits_{S_r(h_0)}F(h)(d\ln f+j\,d\psi)=\lim\limits_{r\to0}
    \oint\limits_{S_r(h_0)}F(h)(d\ln f+j\,d\psi)=F(h_0)\int\limits_{\Psi}(d
    \ln f+j\,d\psi).\]
The integral of the first term cancels due to the bijectivity of
    the function $\ln f$. In this way, we get the following formula
    of the hyperbolic version of the integral Cauchy formula:
\[\oint\limits_{\Gamma}\frac{F(h)}{h-h_0}\,dh=jF(h_0)\int\limits_{\Psi}d\psi.\]
In general, the integral obtained in the right hand side, is
    divergent. However, it can have a meaning, if we introduce the
    formal quantity $\ell_H$ of {\em the size of the hyperbolic space
    of directions} by means of the formula:
\begin{equation}\label{razmer}\frac{\ell_H}{2}\equiv\int\limits_{\mathbb{R}}d\psi.\end{equation}
Taking into consideration the orientations of the pieces $\mathbb{R}$ in
    $\Psi$ (see Fig.\ref{psi}), we get:
\[\int\limits_{\Psi}d\psi=\ell_H/2-\ell_H/2+\ell_H/2-\ell_H/2=0.\]
In this way, the hyperbolic Cauchy formula in some (improper) sense
    has a more simple form than in the complex case:
\begin{equation}\label{coshh1}\oint\limits_{\Gamma}\frac{F(h)}{h-h_0}\, dh=0.\end{equation}

\begin{figure}[bht]\centering\input{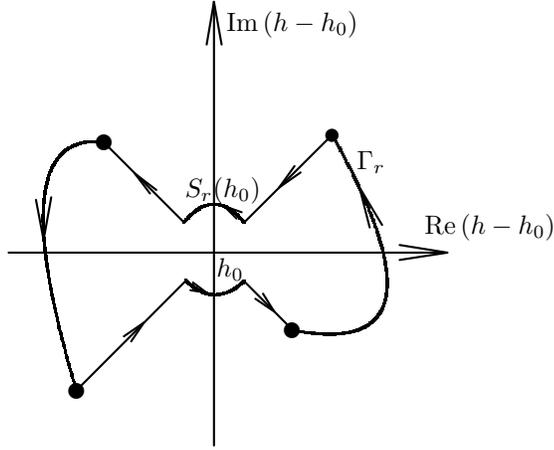}
    \caption{\small Towards the setting the integral Cauchy Theorem
    on the plane of double variable: the contour $\Gamma_r$}.\label{polkont}
\end{figure}

We can obtain a more comprehensive analogue for the standard Cauchy
    formula, if we better examine the closed contour $\Gamma_r$ of
    the form depicted in Fig.\ref{polkont}.
This contour consists of two arcs of arbitrary piecewise-smooth
    simple curves, which lie in the domains $|t-t_0|\ge|x-x_0|$
    and have their ends on the components of the cone $\text{Con}(h_0)$,
    slices of this cone, and a pair of arcs of the Euclidean circle
    of radius $r$ and center at $h_0$, which lie on the components
    of the cone $\text{Con}(h_0)$. The Cauchy-type integral vanishes
    on the contour $\Gamma_r$, with the same general meaning like
    (\ref{coshh1}), due to the fact that the contour $\Gamma_r$
    is homotopic to the initial contour $\Gamma$ inside the domain
    of holomorphy of the function $F(h)/(h-h_0)$. We have now
\begin{equation}\label{coshh2}0=\oint\limits_{\Gamma_r}\frac{F(h)}{h-h_0}\,
    dh=\oint\limits_{S_r(h_0)}\frac{F(h)}{h-h_0}\,
    dh+\oint\limits_{\Gamma'_r}\frac{F(h)}{h-h_0}\, dh,\end{equation}
where $\Gamma'_r\equiv\Gamma_r\setminus S_r(h_0)$. By introducing
    on $S_r(h_0)$ hyperbolic polar coordinates, reiterating the
    previous reasoning, and using the properties of the function
    $f(\psi)$ (its property of being even relative to $\psi$),
    which is provided by the polar equation of the Euclidean circle, we get
\begin{equation}\label{coshh3}
    \lim\limits_{r\to0}\oint\limits_{S_r(h_0)}\frac{F(h)}{h-h_0}\,
    dh=-j\ell_HF(h_0),\end{equation}
whence from (\ref{coshh2}) we infer a more direct analogue of the
    Cauchy formula:
\begin{equation}\label{coshh4}F(h_0)=\frac{1}{\ell_Hj}\oint\limits_{\Gamma_0}
    \frac{F(h)}{h-h_0}\,dh,\end{equation}
where $\Gamma_0=\lim\limits_{r\to0}\Gamma_r$.
    Formally, the obtained formula has a totally equivalent shape
    to the standard complex Cauchy formula, while replacing the
    size of the space of Euclidean directions $\ell_E=2\pi$ with the
    size of the space of hyperbolic directions $\ell_H$
    in the pairs of quadrants whose signs $h\bar h$ coincide.
The quantity $\ell_H$ may be considered as "the fundamental
    constant"\, of the geometry of double numbers.
    While using this constant in computations, we need to
    accurately consider its properties and to use the
    procedure of regularization of expressions.

The hyperbolic Cauchy formula (\ref{coshh4}) can be written in a
    more general form:
\begin{equation}\label{coshhf}F(h_0)=\pm\frac{1}{\ell_Hj}\oint
    \limits_{\Gamma_\pm}\frac{F(h)}{h-h_0}\,dh,\end{equation}
considering the possibility of choosing the contour $\Gamma_-$
    instead of $\Gamma_+=\Gamma_0$ -- obtained by rotating the latter one
    by the Euclidean angle $+\pi/2$, and the reverse orientation of the
    parameter $\psi$ in the domain $|t-t_0|\le |x-x_0|$ relative to the
    general sense of positive (trigonometric, anti-clockwise) tracing of
    the contours in $\H $.

\begin{figure}[bht]\centering\input{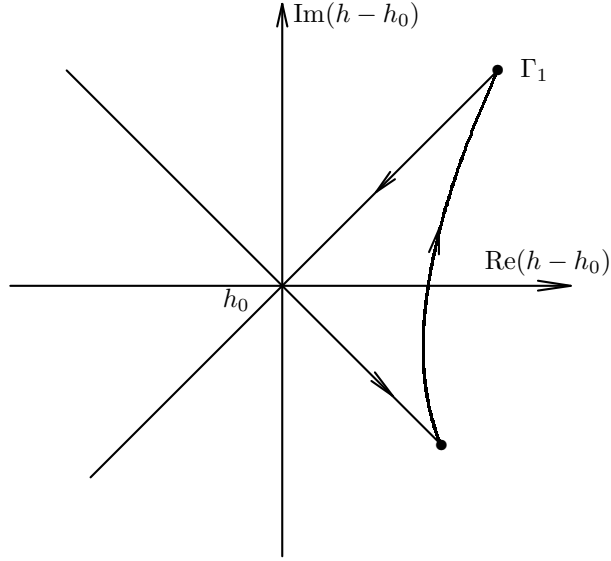}
    \caption{\small Towards the setting the integral Cauchy Theorem
    on the plane of double variable: the contours $\Gamma_n$}.\label{cont3}
\end{figure}

In a similar way one can obtain the following variants of the Cauchy
    formula:
\begin{equation}\label{coshhf1}F(h_0)=(-1)^{n+1}\frac{2}{\ell_Hj}\oint
    \limits_{\Gamma_n}\frac{F(h)}{h-h_0}\,dh,\end{equation}
where $n=1,2,3,4$, the contour $\Gamma_1$ is presented in
    Fig.\ref{cont3}, and the contours $\Gamma_n$ emerge form this
    by rotations of angles $\pi(n-1)/2$ around the point $h_0$.

\bigskip

{\small {\bf Example}.  We illustrate the way in which the Cauchy
    formula works in its form (\ref{coshhf1}) for $n=1$, by means of
    explicitly computing the integral over the contour $\Gamma_1$. We
    obtain
\begin{equation}\label{ex1}\frac{2}{\ell_Hj}\oint\limits_{\Gamma_1}\frac{F(h)}{h-h_0}\,dh
    =\frac{2}{\ell_Hj}\oint\limits_{\Gamma_1}\frac{F(h)-F(h_0)}{h-h_0}\,
    dh+\frac{2F(h_0)}{\ell_Hj}\oint\limits_{\Gamma_1}\frac{dh}{h-h_0}.
\end{equation}
%

The expression within the first integral is a holomorphic mapping in
    the domain which is bounded by the contour $\Gamma_1$ and on the
    contour itself, and therefore this integral cancels. On the cone
    $\text{Con}_+^\uparrow$ we choose as integration variable
    $t\in[t_0+\tau_1,t_0]$, and on the cone $\text{Con}_-^\uparrow$,
    $t\in[t_0,t_0+\tau_2]$, where $\tau_1$ and $\tau_2$ are the
    abscissae of the endpoints of the curvilinear part of the contour
    $\Gamma_1$ (the upper and lower ends, respectively) in the system
    of coordinates having the origin at the point $h_0$. In this way,
    the main contribution of the integral is:
\[\frac{2F(h_0)}{\ell_Hj}\left[\int\limits_{t_0+\tau_1}^{t_0}
    \frac{(1+j)\,dt}{(1+j)(t-t_0)}+\int\limits_{t_0}^{t_0+\tau_2}
    \frac{(1-j)\,dt}{(1-j)(t-t_0)}\right]=\]
\[\frac{2F(h_0)}{\ell_Hj}(\ln 0-\ln\tau_1+\ln\tau_2-\ln 0)
    =\frac{2F(h_0)}{\ell_Hj}\ln(\tau_2/\tau_1)=0.\]
In the prior to last equality we took into consideration
    the pairwise canceling of two logarithmic singular terms,
    and in the last one, we used the "unboundedness property"
    of the fundamental constant $\ell_H$. In this way,
    the main contribution in the Cauchy integral is provided only by the segment
    $\Gamma'$ of the contour located between the components of the
    cone $\text{Con}(h_0)$. Passing to polar coordinate system with center
    at the point $h_0$, we get:
\[\frac{2F(h_0)}{\ell_Hj}\int\limits_{\Gamma'}\frac{dh}{h-h_0}
    =\frac{2F(h_0)}{\ell_Hj}\int\limits_{\Gamma'}(d\ln\varrho+d\psi).\]
The integral in the first term cancels, due to the fact that at the
    ends of the contour $\Gamma'$ we have $\varrho=0$. We integrate
    the second term by considering (\ref{razmer}), and we get $F(h_0)$,
    fact which confirms the validity of the hyperbolic Cauchy
    formula in the form (\ref{coshhf1}).}\par\bigskip

\begin{figure}[htb]\centering\input{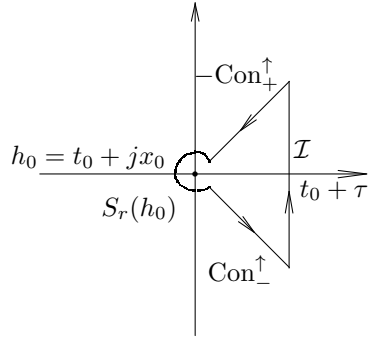}
\caption{\small Towards the question of finding the coefficients
    of the Taylor series for $h$-holomorphic mappings. }\label{g1}
\end{figure}

We shall examine now the question of possibility of computing the
    coefficients of the Taylor series of a $h$-holomorphic function,
    by means of a formula similar to Cauchy one in complex analysis.
To this aim, we consider the integral of the form:
\begin{equation}\label{idh}\oint\limits_{\Gamma}(h-h_0)^\alpha\, dh,\quad \alpha\in R.\end{equation}
We deform the contour $\Gamma$ in such a way, that it gets the form
    $\Gamma'$, shown in Fig.\ref{g1} (the value of the integral
    remaining, obviously, unchanged).

Splitting the integral into additive terms which correspond to
    different segments of the contour $\Gamma'=-\text{Con}_+^\uparrow(h_0)\cup
    S_r(h_0)\cup\text{Con}_-^\uparrow(h_0)\cup\mathcal{I}$, we get
    for contributions over cones (the notations being similar to the ones
    used in the previous example, with $\tau_1=\tau_2$):
\[\int\limits_{-\text{Con}_+^\uparrow(h_0)}\hspace{-1em}(h-h_0)^\alpha\,
    dh+\int\limits_{\text{Con}_-^\uparrow(h_0)}\hspace{-1em}(h-h_0)^\alpha\,dh
    =\int\limits_{t_0+\tau}^{t_0+r}(t-t_0)^\alpha(1+j)^{\alpha+1}\,dt+
    \int\limits_{t_0+r}^{t_0+\tau}(t-t_0)^\alpha(1-j)^{\alpha+1}\,dt.\]
Using the identity (\ref{eqqh1}) and collecting the similar terms,
    we obtain the following result:
\begin{equation}\label{conr}\int\limits_{\text{Con}_-^\uparrow(h_0)-
    \text{Con}_+^\uparrow(h_0)}\hspace{-2em}(h-h_0)^\alpha\,dh
    =j\frac{2^{\alpha+1}}{\alpha+1}[r^{\alpha+1}-\tau^{\alpha+1}].\end{equation}
The integral over the straight line segment is easy to compute as
    well:
\begin{equation}\label{otrez}\int\limits_{\mathcal{I}}(h-h_0)^\alpha\,dh|_{h=h_0+\tau+jsx}
    =\int\limits_{-1}^1(\tau+js\tau)^{\alpha}j\tau\,ds
    =\tau^{\alpha+1}\int\limits_{1-j}^{1+j}\xi^\alpha\,d\xi
    =j\frac{2^{\alpha+1}\tau^{\alpha+1}}{\alpha+1}.\end{equation}
Comparing (\ref{conr}) with (\ref{otrez}), we conclude that
\[\int\limits_{\Gamma'}(h-h_0)^\alpha\,dh
    =\int\limits_{S_r(h_0)}\hspace{-1em}(h-h_0)^\alpha\,
    dh+j\frac{2^{\alpha+1}}{\alpha+1}r^{\alpha+1},\ (\alpha\neq-1).\]
For the integral over the circle $S_r(h_0)$, by means of the
    parametrization (\ref{regih}), we obtain the representation:
\[\int\limits_{S_r(h_0)}\hspace{-1em}(h-h_0)^\alpha\,dh
    =[(-1)^{\alpha+1}-j^{\alpha+1}-(-j)^{\alpha+1}]\int
    \limits_{-\infty}^{+\infty}\varrho^{\alpha+1}
    e^{j(\alpha+1)\psi}(d\,\ln\varrho+j\, d\psi).\]
By substituting $\varrho=\varrho(\psi)=rf(\psi)$ (for us, the
    concrete form of $f$ is of less importance, and what matters is the
    dependence of this mapping only on $\psi$, and its independence of
    $r$), we get the expression:
\[\int\limits_{S_r(h_0)}\hspace{-1em}(h-h_0)^\alpha\,dh=r^{\alpha+1}K_1(\alpha)\]
and hence the overall result:
\begin{equation}\label{idhh}\oint\limits_{\Gamma_1}(h-h_0)^\alpha\,dh
    =r^{\alpha+1}K_2(\alpha),\ (\alpha\neq-1).\end{equation}
where $K_1(\alpha)$ and $K_2(\alpha)$ are certain functions, which
    depend only on $\alpha$. Due to the homotopicity of the
    contours having different values for $r$, the expression
    (\ref{idhh}) should not depend on $r$. This is possible only if
    $K_2(\alpha)=0$ for $\alpha\neq-1$. Taking into consideration
    the previous result for $\alpha=-1$, we infer the following final formula:
\begin{equation}\label{rezt}\oint\limits_{\Gamma}(h-h_0)^\alpha\, dh
    = \left\{\begin{array}{lr}0,&\alpha\neq-1;\\j\ell_H,&\alpha=-1.\end{array}\right.
\end{equation}

%
%
\subsection{Conformal-analytic hyperbolic mappings}

Our formerly introduced bilinear form $\mathcal{G}$ behaves under
    $h$-holomorphic mappings $F(h)$ as a relative scalar:
\begin{equation}\label{Bh}\mathcal{G}\mapsto \mathcal{G}'=|F'(h)|^2\mathcal{G},\end{equation}
where $F'(h)=dF/dh$, whence there follow the transformation laws
    for the hyperbolic elements of length and area:
\begin{equation}\label{confh}dl'=|F'|\, dl;\quad (dh\wedge d\bar h)'
    =|F'|^2(dh\wedge d\bar h).\end{equation}
Like in the complex case, the property of being a relative scalar of
    the area form takes place for any diffeomorphism, and the first
    equality in (\ref{confh}) means, that the $h$-holomorphic functions
    form the conformal mappings of the hyperbolic plane, i.e., they
    preserve the hyperbolic angles between curves at each point, where $|F'(h)|^2\neq0.$
This fact is tightly related to the already proved invariance of
    cones $\text{Con}$ relative to $h$-holomorphic mappings.
We note, that $|F'|^2=|\bigtriangledown u|^2=|\bigtriangledown v|^2
    =\Delta_F$, where $\bigtriangledown$ is the gradient
    operator for the pseudo-Euclidean metric, and $\Delta_F$ is the
    Jacobian of the mapping $F$, regarded as mapping $R^2\to R^2$.

Like in the complex case, each diffeomorphism $f:\R^2\to\R^2$ can be
    regarded as a smooth vector field. The vector fields which
    correspond to $h$-holomorphic mappings of double variable have
    certain interesting and important applicative properties.
    From the condition (\ref{crh}) it follows, that each component
    of the vector field $F(h)=U+iV$ is a hyperbolic $h$-harmonic function,
    i.e., it satisfies the wave-equation of second order (\ref{harmh}).
    The $h$-harmonic functions which are mutually related by the
    Cauchy-Riemann conditions (\ref{crh}), will be called
    {\em $h$-conjugate}. Any $h$-harmonic mapping on the Cartesian plane
    defines (up to a constant) its hyperbolic conjugate partner.
    The hyperbolic Cauchy-Riemann conditions (\ref{crh}) have, in
    terms of vector analysis in hyperbolic space, the following
    geometric meaning: {\em the vector field $\bar F=U-jV$ is $h$-potential
    and $h$-solenoidal}, i.e., the components $\{U,-V\}$ of the vector field
    $F$ satisfy the relations:
\begin{equation}\label{vecth}\bar F_{,h}=0\Leftrightarrow\text{roth}\,
    F\equiv U_{,x}-V_{,t}=0;\quad \text{divh}\,F\equiv U_{,t}-V_{,x}=0.\end{equation}
The physical meaning of these conditions and the corresponding
    problems with initial/boundary conditions, which are naturally
    solved using hyperbolic conformal mappings, will be further
    considered in a separate section. We note here, that the family
    of curves $U=\text{const}$ and $V=\text{const}$ define on the
    Cartesian plane $\R^2$ pseudo-orthogonal families of curves,
    where $\bigtriangledown U\cdot\bigtriangledown V\equiv
    U_{,t}V_{,t}-U_{,x}V_{,x}=0$ holds everywhere true, due to
    the conditions (\ref{crh}).
%
%
\section{Properties of some elementary functions of double \\variable}\label{elemmm}

We examine now in detail the properties of the basic elementary functions of double variable.
%
%
\subsection{Power functions $F(h)=h^n$}

Unlike the power function of complex variable, the cases of even powers $n$ and odd powers $n$
    are essentially different.
Indeed, when passing to the exponential representation (\ref{exph}), we get:
\begin{equation}\label{step1}h=\epsilon\varrho e^{j\psi}\mapsto \epsilon^n\varrho e^{jn\psi}\end{equation}
Since for any even $n$ we have $\epsilon^n=1$, we conclude that {\em the power mapping
    $h\mapsto h^n$, for $n=2k$, $k\in\Z$ bijectively maps each of the quadrants I, II, III and
    IV on the first quadrant, and maps the cones $\text{Con}_\pm\to\text{Con}_{\pm}$}.
    On the contrary, {\em for odd $n$, each of the coordinate quadrants is bijectively mapped
    by the transformation $h\mapsto h^n$ $n=2k+1$, $k\in\Z$ into itself}.
We easily see from (\ref{step1}) that the net of coordinate lines
    $\varrho=\text{const}$, $\psi=\text{const}$ maps to the net of coordinate lines
    $\varrho'=\varrho^n=\text{const}$, $\psi'=n\psi=\text{const}$ for all integer $n$.
In the case of positive integers $n$, the radial lines extend for $\varrho>1$
    and shrink for $\varrho<1$. Moreover, they rotate from the value $\psi=0$ in the sense
    of the cone components which correspond to their signs.
For integer negative $n$, there occur a complimentary inversion relative to the unit spheres
    $\varrho=1$ and a mirroring of the space of angles $\Psi\to -\Psi$. As example of function
    with even $n$ we examine the mapping $w=h^2=x^2+y^2+2jxy=\varrho^2e^{2\psi}$.

\begin{figure}[htb]\centering\input{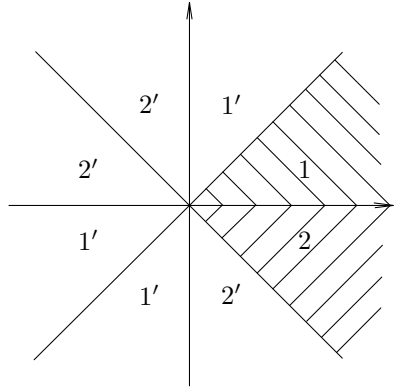}
    \caption{\small The global structure of the mapping $h\mapsto h^2$}.\label{quad}
\end{figure}

In Fig.\ref{quad} there is represented the global structure of the
mapping $h\mapsto h^2$:
    the quadrant 1-2 is mapped into itself (its borders map to the corresponding ones),
    and the mapping of the other quadrants into quadrant 1-2 is shown by the corresponding
    figures (accented figures which identify a quadrant, show how the corresponding quadrant
    maps to the quadrant 1-2). In this way, the mapping $h\mapsto h^2$ is 4-fold.
A similar situation occurs with the mapping: $h\to h^{2k}$ $k\in\Z$.

We present in Fig.\ref{mapindex1} an illustrative representation of mapping $h\to h^2$.

{\centering\small\refstepcounter{figure}\label{mapindex1}
    \includegraphics[width=.3\textwidth]{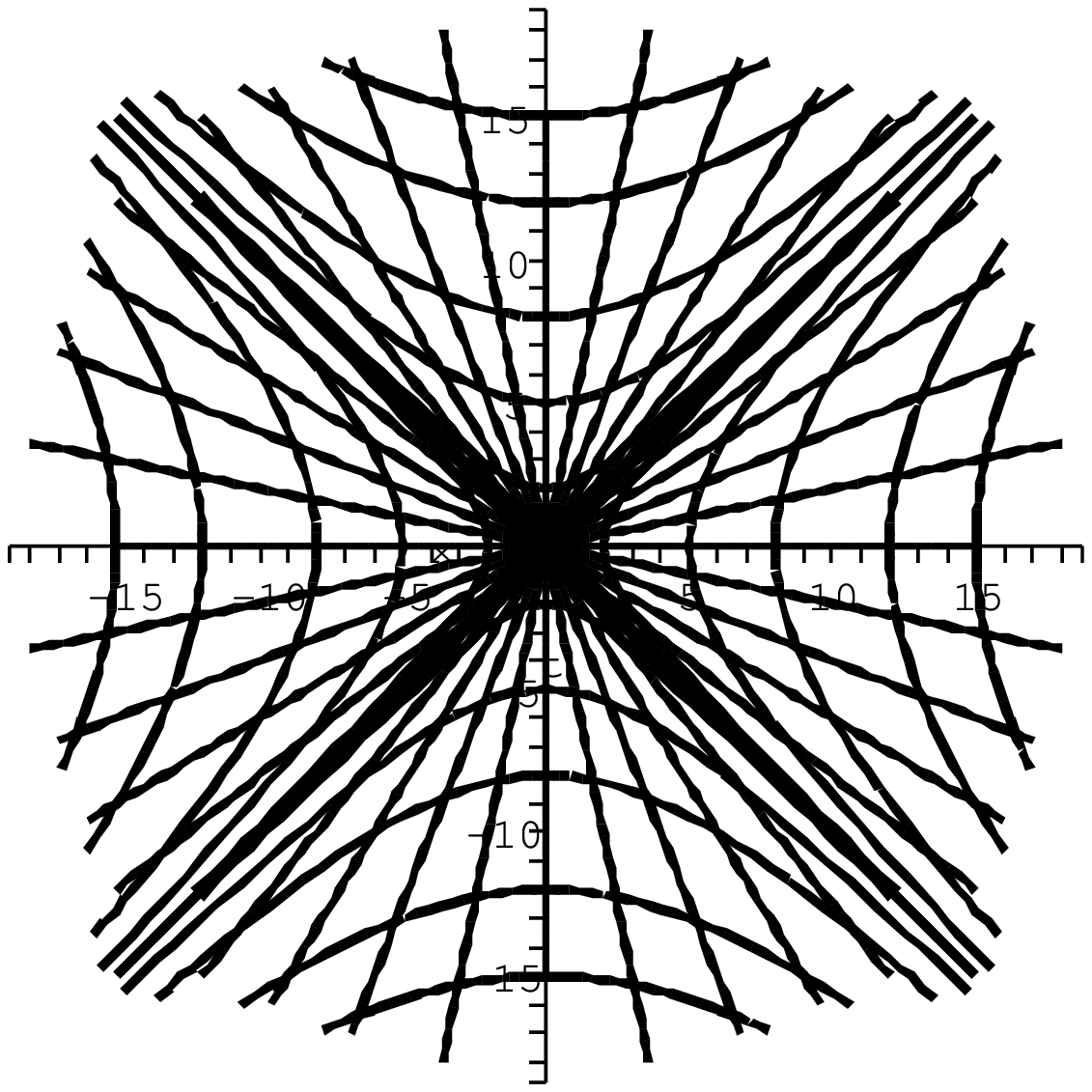}\includegraphics[width=.3\textwidth]{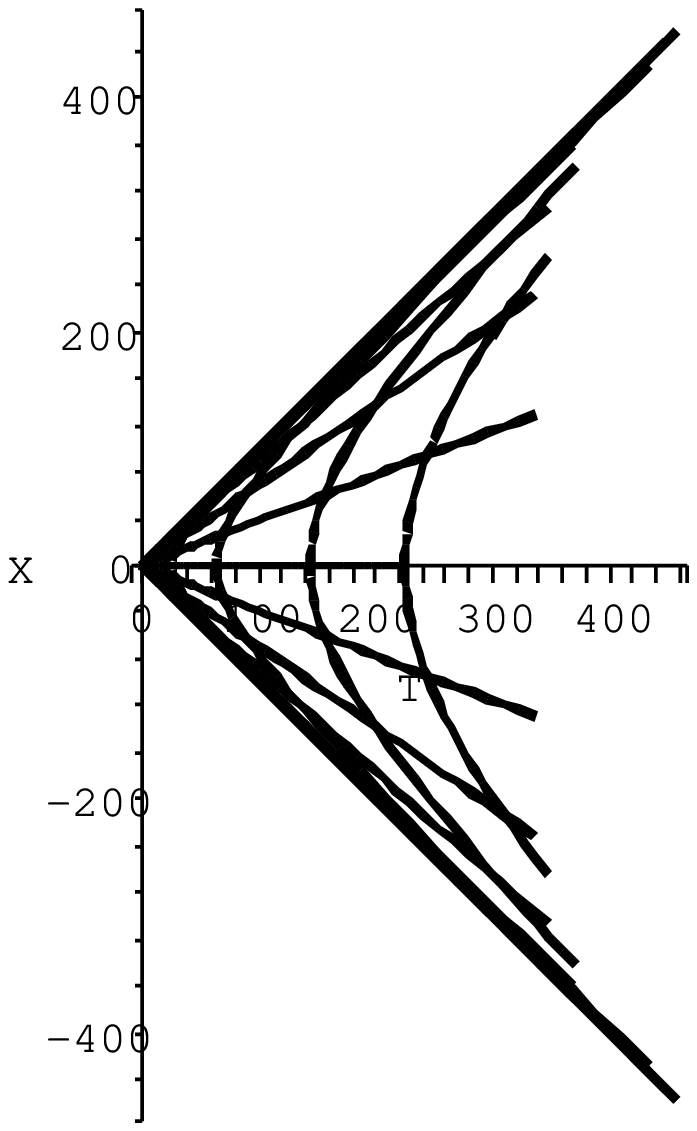}
    \medskip\nopagebreak\\ Fig.~\thefigure. The hyperbolic polar system of coordinates (left)
    and the image of its first quadrant under the mapping $h\mapsto h^2$. }


\par\medskip

\begin{figure}[htb]\centering\input{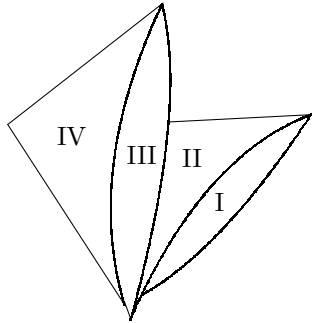}
    \caption{\small The hyperbolic Riemannian surface of the 4-valued mapping $h\mapsto h^{1/2k}$,
    $k\in\Z$}.\label{riemind}
\end{figure}

From the properties of power functions it is easy to derive the properties of roots of different
    orders and of the rational powers $h\mapsto h^{1/n}$ and $h\mapsto h^{m/n}$.
Each root $\sqrt[n]{h}$ of even order is defined in quadrant I. Such a root is a 4-valued
    function.
Each leaf of the hyperbolic Riemannian surface of this function represents a unit copy of the
    first quadrant I, as shown on Fig.\ref{quad}. On each leaf, the mapping is bijective.
All the leaves glue together into a Riemannian surface, which represents $\R^2$,
    with the point $(0;0)$ belonging to all the leaves and being a hyperbolic analogue of
    the branching point. The Riemannian surface of the roots of even orders can be illustrated
    by means of a sheet of paper, folded as shown on Fig.\ref{riemind}.

The roots of odd order are bijective in each of the 4 quadrants.
%
%
\subsection{The exponential of double variable $w=e^h$}\label{exppp}

The relations $e^h=e^{t+jx}=e^te^{jx}$ naturally lead to the global structure of the
    exponential mapping, which is represented in Fig.\ref{exp}.

{\centering\small\refstepcounter{figure}\label{exp}
    \includegraphics[width=.4\textwidth]{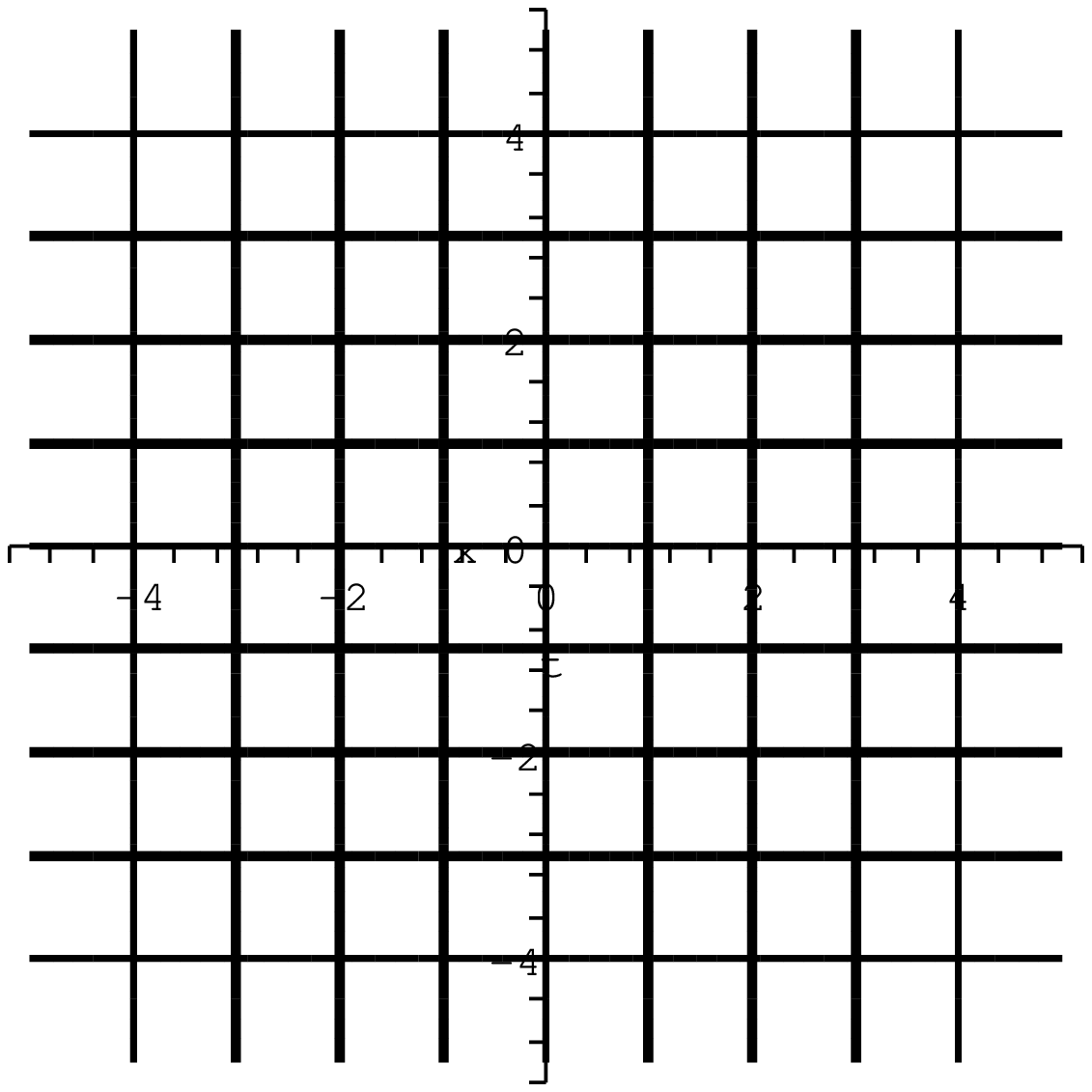}\includegraphics[width=.4\textwidth]{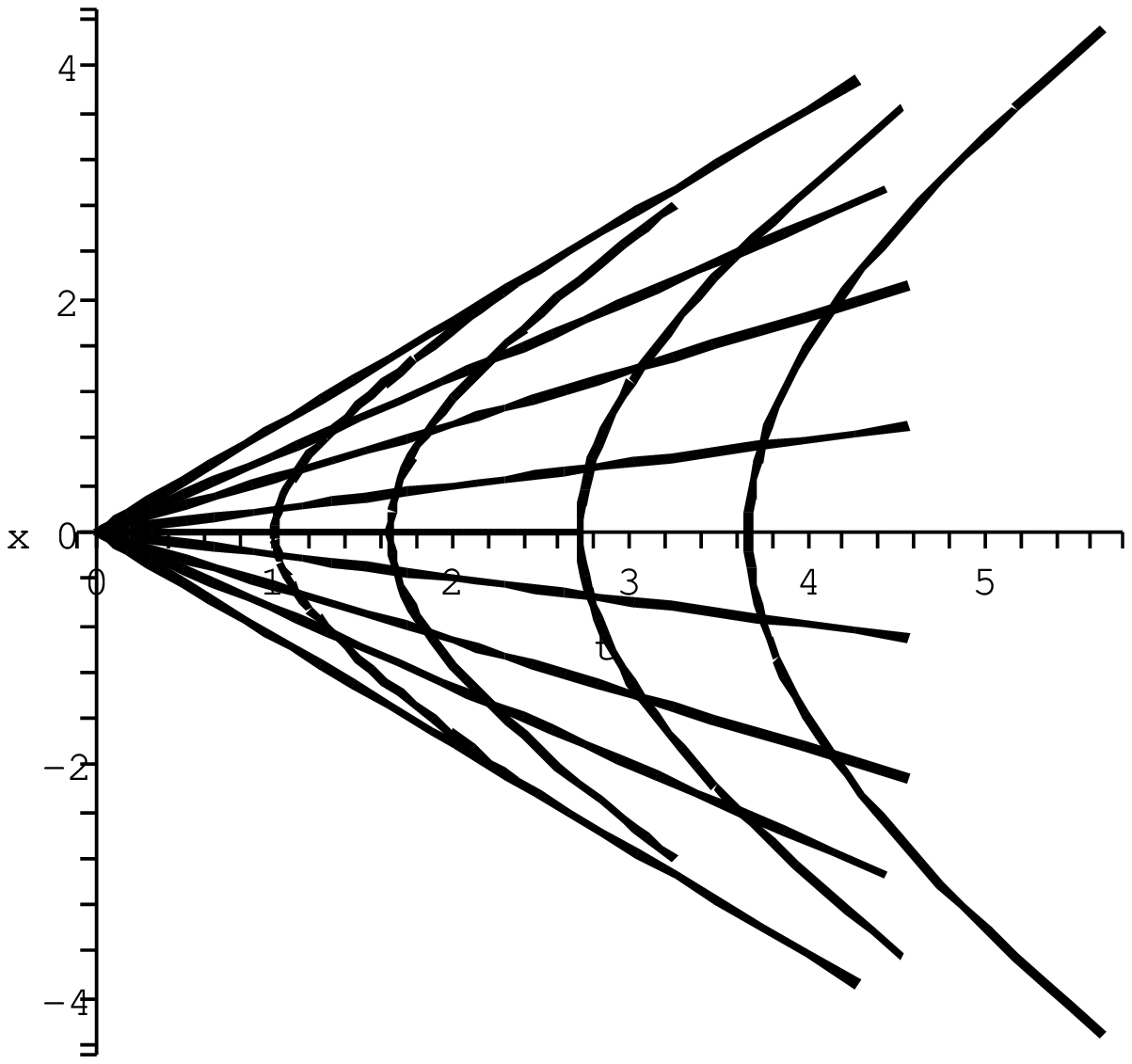}
    \medskip\nopagebreak\\ Fig.~\thefigure. The global structure of the mapping $h\mapsto e^h$}.\par\medskip

The rectangular pseudo-orthogonal net on the plane with variable $h$ is mapped by the
    exponential into the pseudo-orthogonal net, which consists of the rays and hyperbolas
    from the first quadrant and the vertex at the point $h=0$. The mapping $h\mapsto e^h$ is bijective.
Obviously, the inverse mapping $\ln h=\ln\varrho+j\psi$ is defined inside the first quadrant.
    On its border (i.e., on the cone $\text{Con}^{\uparrow}(0)$) the polar coordinate system is
    ineffective and we need supplementary investigation regarding the behavior of the mapping
    $h\mapsto e^h$, which we shall not address here.
%
%
\subsection{The trigonometric functions $\sin h$, $\cos h$ and their inverses}

While writing the sine of double variable:
\begin{equation}\label{sinn}\sin h=\sin(t+jx)=\sin t\cos x+j\sin x\cos t,\end{equation}
we note that the lines $x=\text{const}$ and $t=\text{const}$ mapped into families of
    ellipses with their center at the point $(0;0)$. These lines are wrapped around
    these ellipses countless times. In Fig.\ref{sin} (at right) there are shown the
    images of the squares with various sides, having the center at the point $(0;0)$ (at left).

{\centering\small\refstepcounter{figure}\label{sin}
    \includegraphics[width=.4\textwidth]{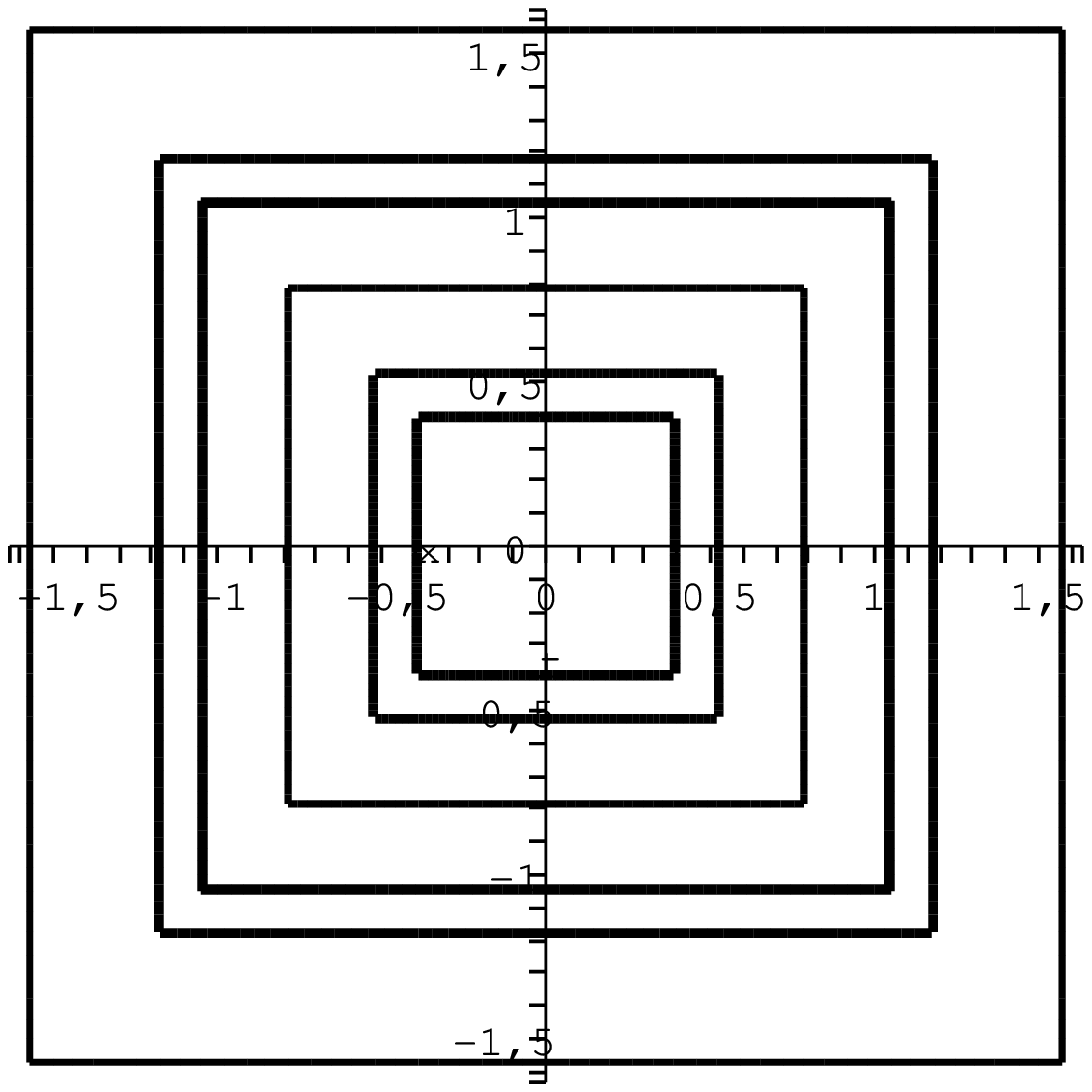}\includegraphics[width=.4\textwidth]{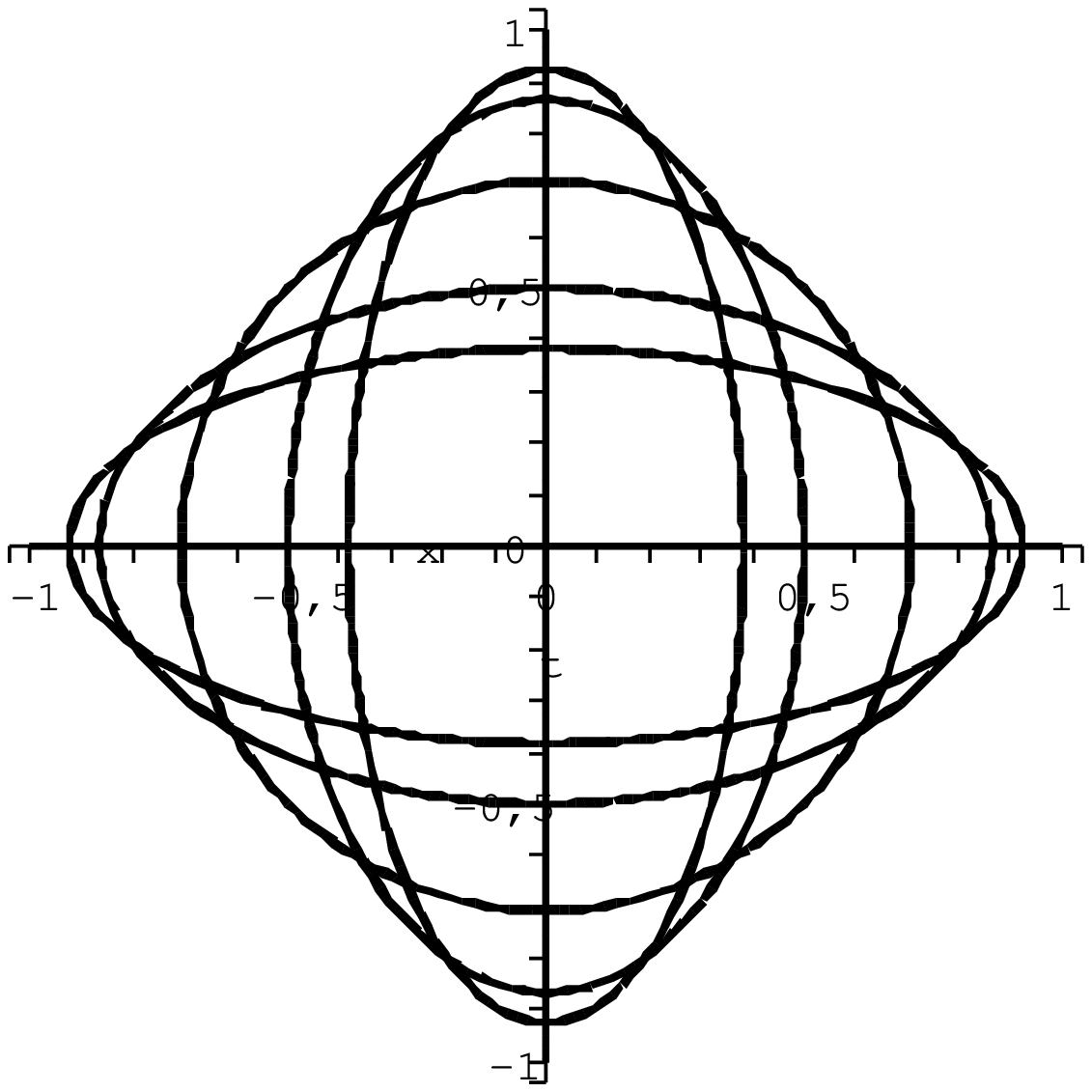}
    \medskip\nopagebreak\\ Fig.~\thefigure. To the properties of the mapping $h\mapsto\sin h$. }\par\medskip

Each square is mapped into a 4-ray star-shaped figure, and the
square with the side $\pi/2$
    is mapped into a circle, the square with the side $\pi$ is mapped  into a coordinate cross
    with the vertices at the points $(1;0)$, $(0;1)$, $(-1;0)$, $(0;-1)$.

In order to clear global structure of the mapping $h\mapsto \sin
h$ it is more convenient to consider a system of fundamental
squares. One of them is presented in Fig. \ref{sina} (from the
left).

{\centering\leftskip0em\rightskip5em\small
\refstepcounter{figure}\label{sina}
\includegraphics[width=.4\textwidth]{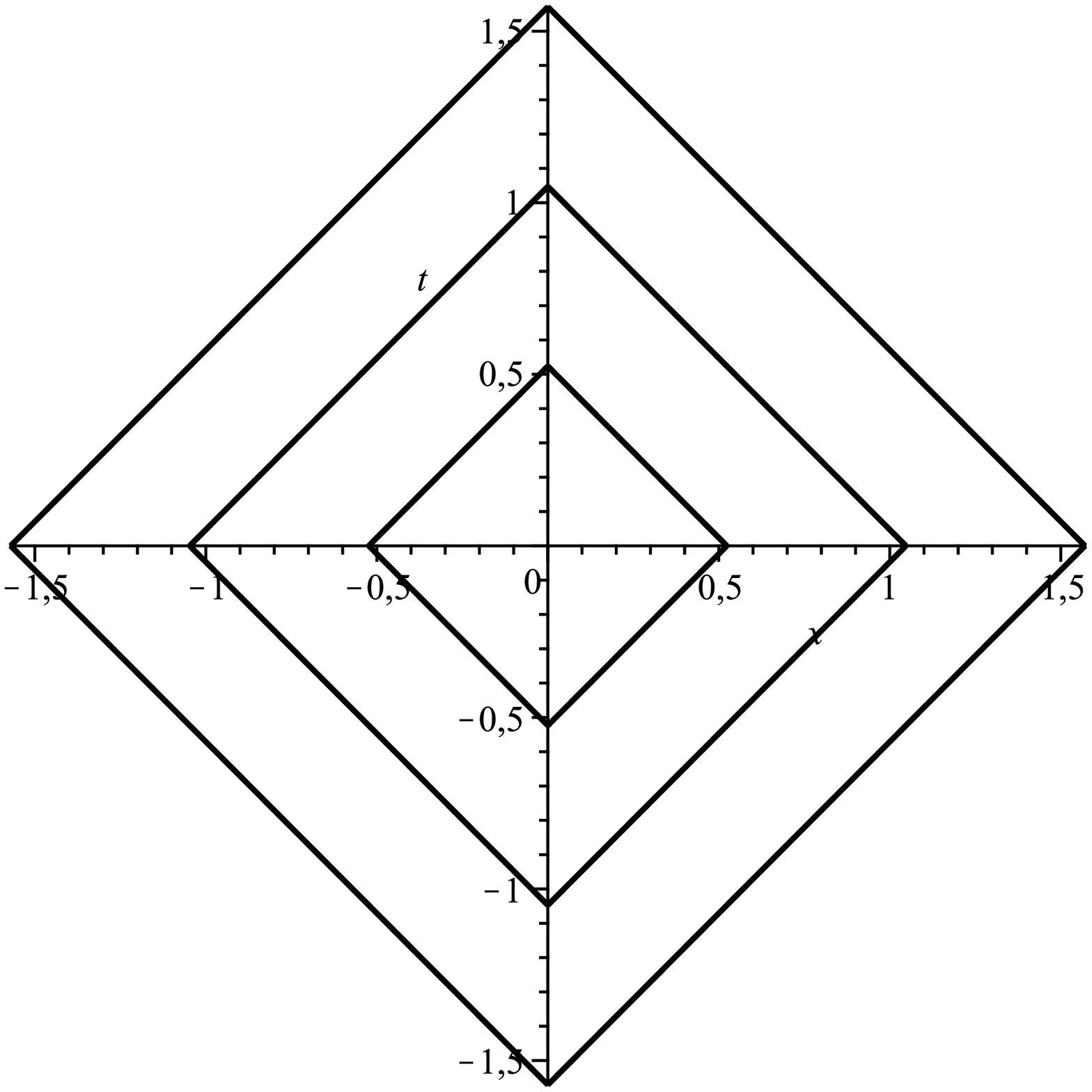}
\includegraphics[width=.4\textwidth]{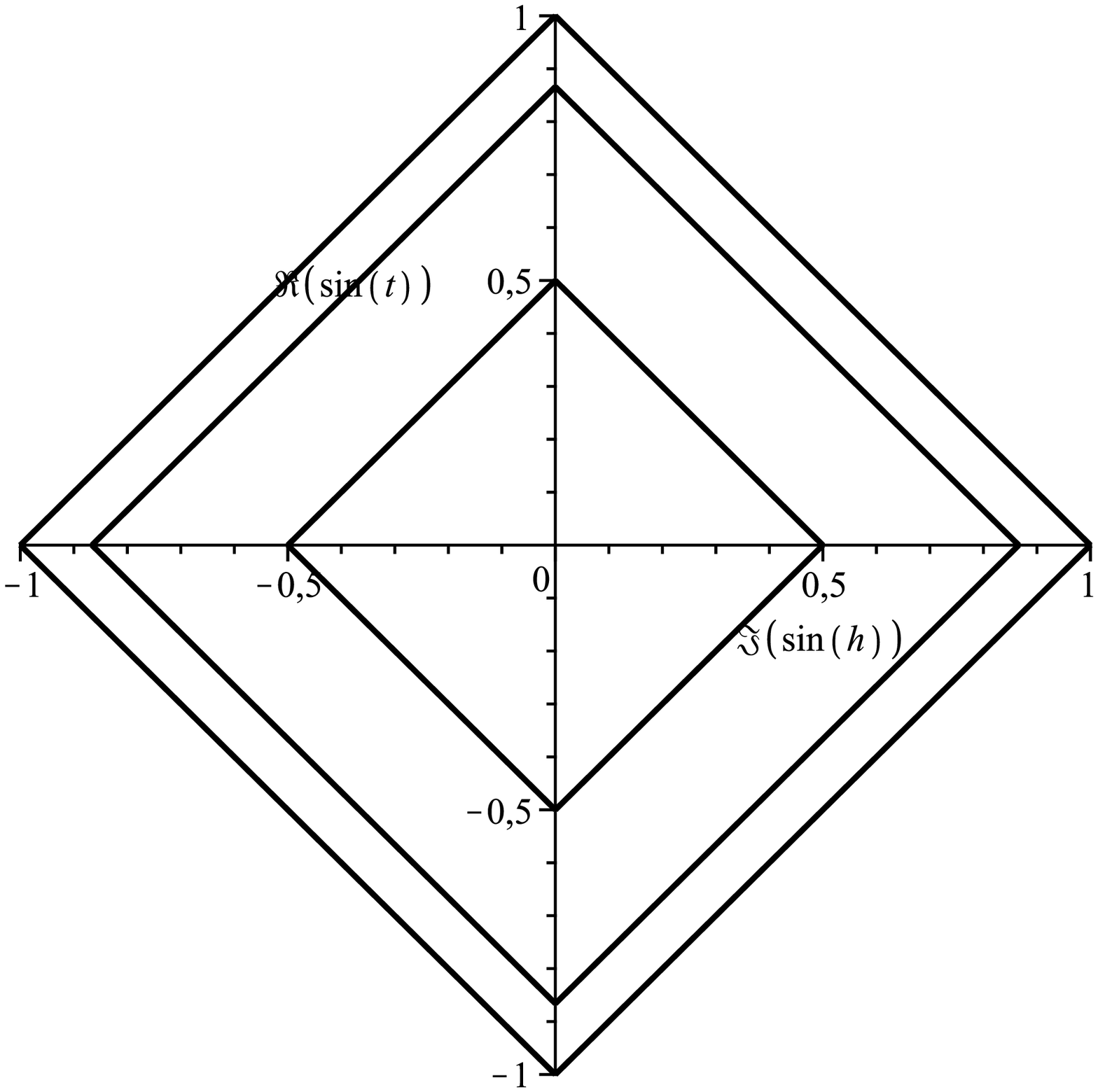}
\medskip
\nopagebreak\par Fig.~\thefigure. Global structure of the mapping
$h\mapsto \sin h.$ Fundamental square (the largest from the left)
is mapped onto the largest square from the right. In each
fundamental square function $\sin h$ is bijective.}

\bigskip

The whole plane of variable $h$ is covered by the squares such as
in fig. \ref{sina} (from the left) by the shifts on vectors $(k\pm
jm)\pi/2,$ $k,m\in \mathbb{Z}.$ Note, that the mapping
$h\mapsto\sin h$ in two neighboring squares have opposite screw
sense (i.e. opposite sign of Jacobian).
%

It is easy to see, that the mapping $h\mapsto\cos h$ similarly works, with the difference that
    the whole family of "fundamental squares"\, is left-shifted on the plane of variable $h$ with
    $\pi/2$ (since $\cos h=\sin(h+\pi/2)$.)

So, the function $\arcsin$ (and $\arccos$) we can define on the
square
    with vertices at the point $(1;0),(0;1), (-1;0),(0;-1)$ (on such a square, left-shifted by
    $\pi/2$ ). Explicit formulas for arcsin and arccos have the form:
\[\begin{array}{ll}\arcsin\,h=&\frac{1}{2}[\arcsin((t+x)\sqrt{1-(t-x)^2}+(t-x)\sqrt{1-(t+x)^2})\smallskip\\
        &+j\arcsin((t+x)\sqrt{1-(t-x)^2}-(t-x)\sqrt{1-(t+x)^2})];\medskip\\
    \arccos\,h=&\frac{1}{2}[\arccos(t^2-x^2-\sqrt{1-(t-x)^2}\sqrt{1-(t+x)^2})\smallskip\\
        &+j\arccos(t^2-x^2+\sqrt{1-(t-x)^2}\sqrt{1-(t+x)^2})].\end{array}\]
%
%

%

\medskip

%
%
\subsection{The trigonometric functions $\tan h$, $\cot h$ and its inverses}

By identifying in the functions $w=\tan h$ the real and the imaginary parts, after
    elementary transformations, we get:
\[\tan h=\frac{\sin 2t+ j\sin 2x}{\cos 2t+\cos 2x}.\]
This mapping transforms the square with the center in the point
$(0;0)$ and edge $\pi/2$ into the
    domain bounded by hyperbolas, and the rectangular net from the original square, into the
    symmetric hyperbolic net inside the domain (Fig.\ref{tan}).

{\centering\small\refstepcounter{figure}\label{tan}
    \includegraphics[width=.4\textwidth]{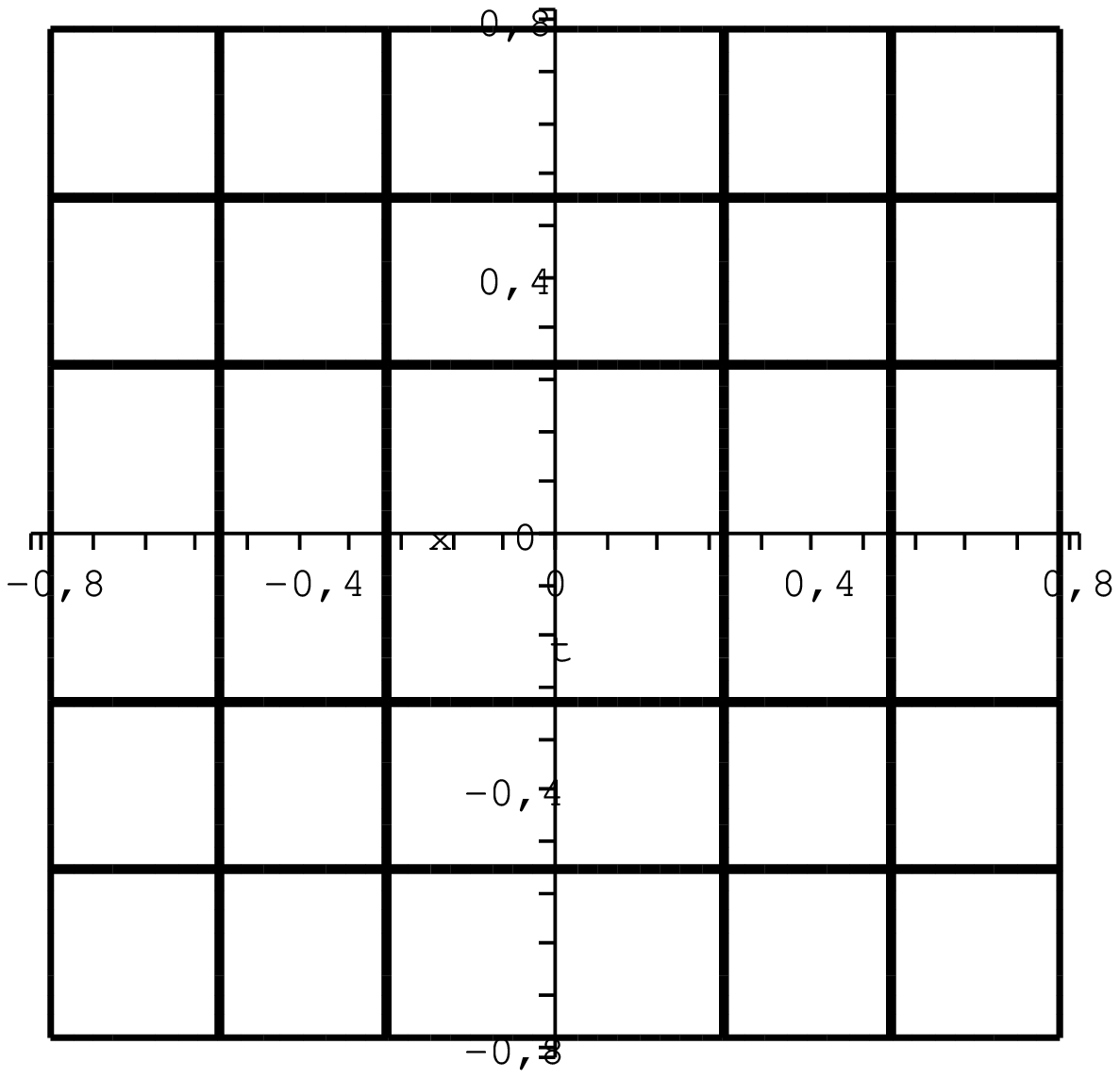}\includegraphics[width=.4\textwidth]{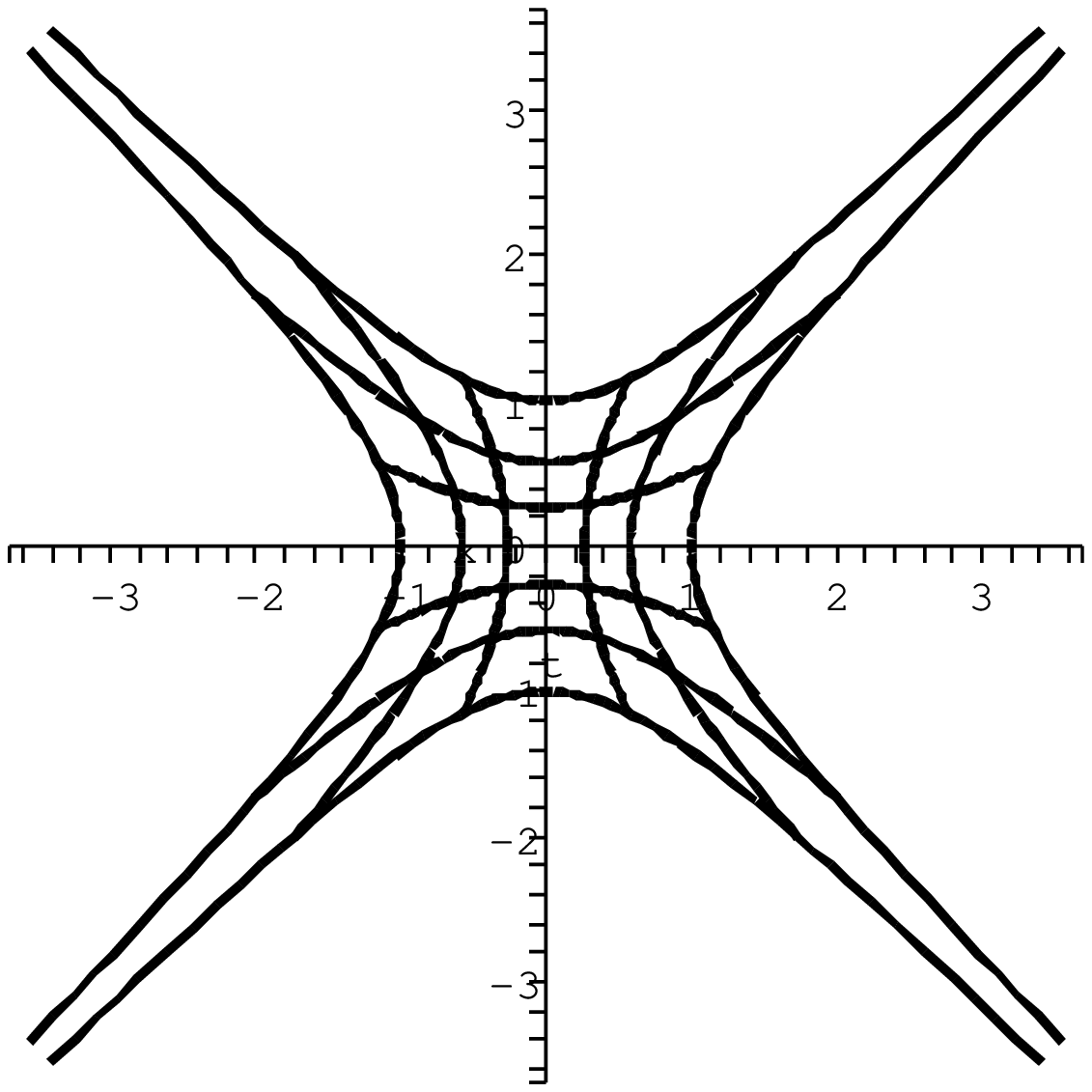}
    \medskip\nopagebreak\\ Fig.~\thefigure. The structure of the mapping $h\mapsto \tan h$}.\par\medskip

Globally, the mapping $h\to \tan h$ has an infinite number of leaves. These are the squares which are
    obtained from the fundamental square $(\pi/2;0),(0;\pi/2),(-\pi/2;0),(0;-\pi/2)$ (see left fig. \ref{sina}) by means of
    translations with vectors multiple of $\pi$ by $t$ and by $x$. 
%

In view of the identity $\cot h=-\tan(h-\pi/2)$, in a similar way behaves the mapping $w=\cot h$.
    The functions $\arctan$ and $\text{arccot}$ are multi-valued, their independent branches
    can be identified in each of the fundamental squares. E.g., the mapping $\arctan\, h$ has
    in coordinates the following explicit form:
\[\arctan\,h=\frac{1}{2}\left\{\arctan\left[\frac{2t}{1-t^2+x^2}\right]+
    j\arctan\left[\frac{2x}{1+t^2-x^2}\right]\right\}.\]
%
%
%
\subsection{The hyperbolic functions $\sinh h$, $\cosh h$, $\tanh h$, $\coth h$ and their inverses}

Separating, like in the case of the elliptic sine, in the function $w=\sinh h$
    the real and the imaginary parts, we obtain the expression:
\[\sinh h=\sinh t\cosh x+j\sinh x\cosh t.\]
%
%

It is easy to see that the rectangular coordinate net $(t,x)$ is mapped into the hyperbolic set
    on the plane of images $w$ (Fig.\ref{hyps}).

{\centering\small\refstepcounter{figure}\label{hyps}
    \includegraphics[width=.3\textwidth]{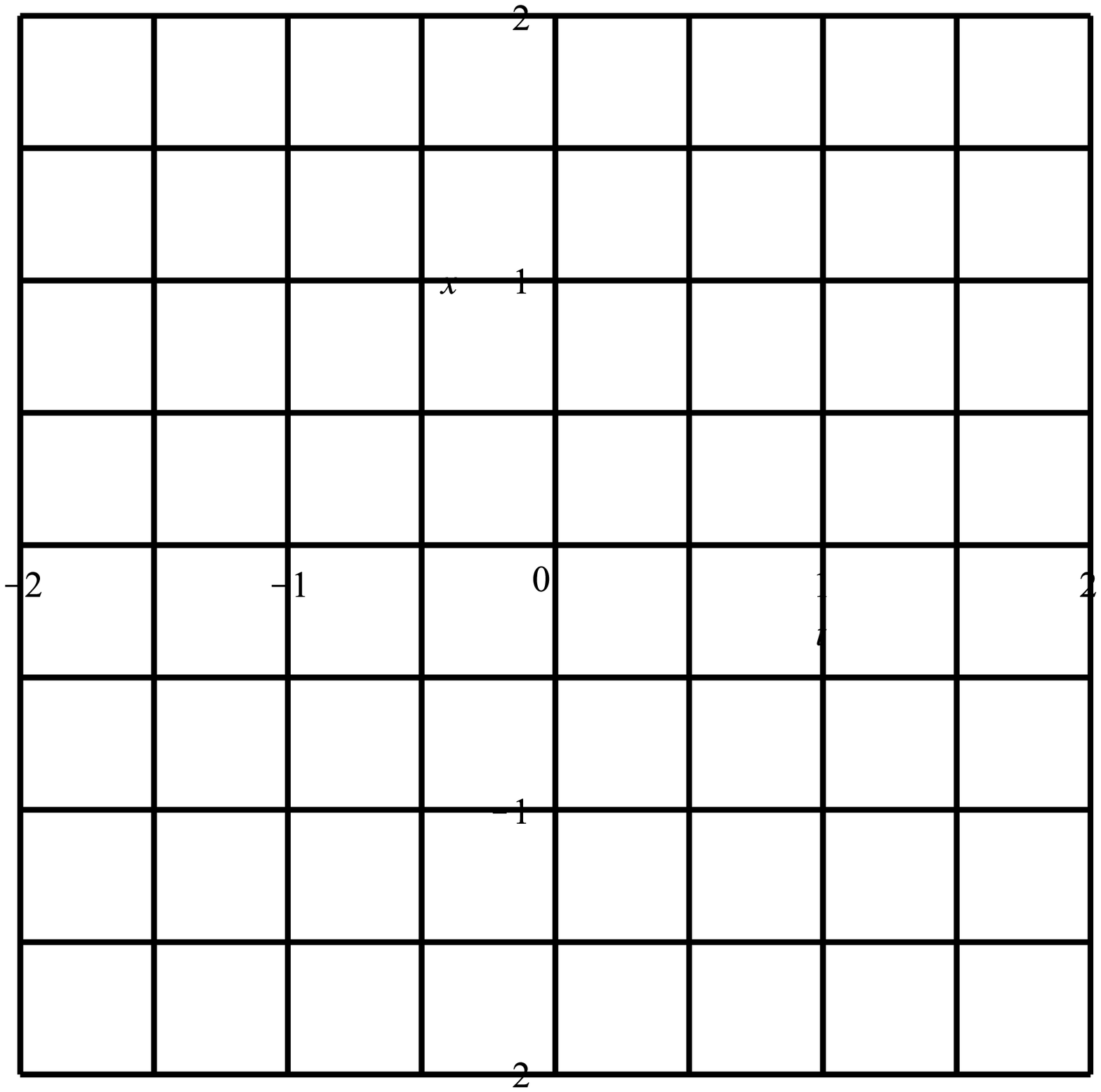}\includegraphics[width=.3\textwidth]{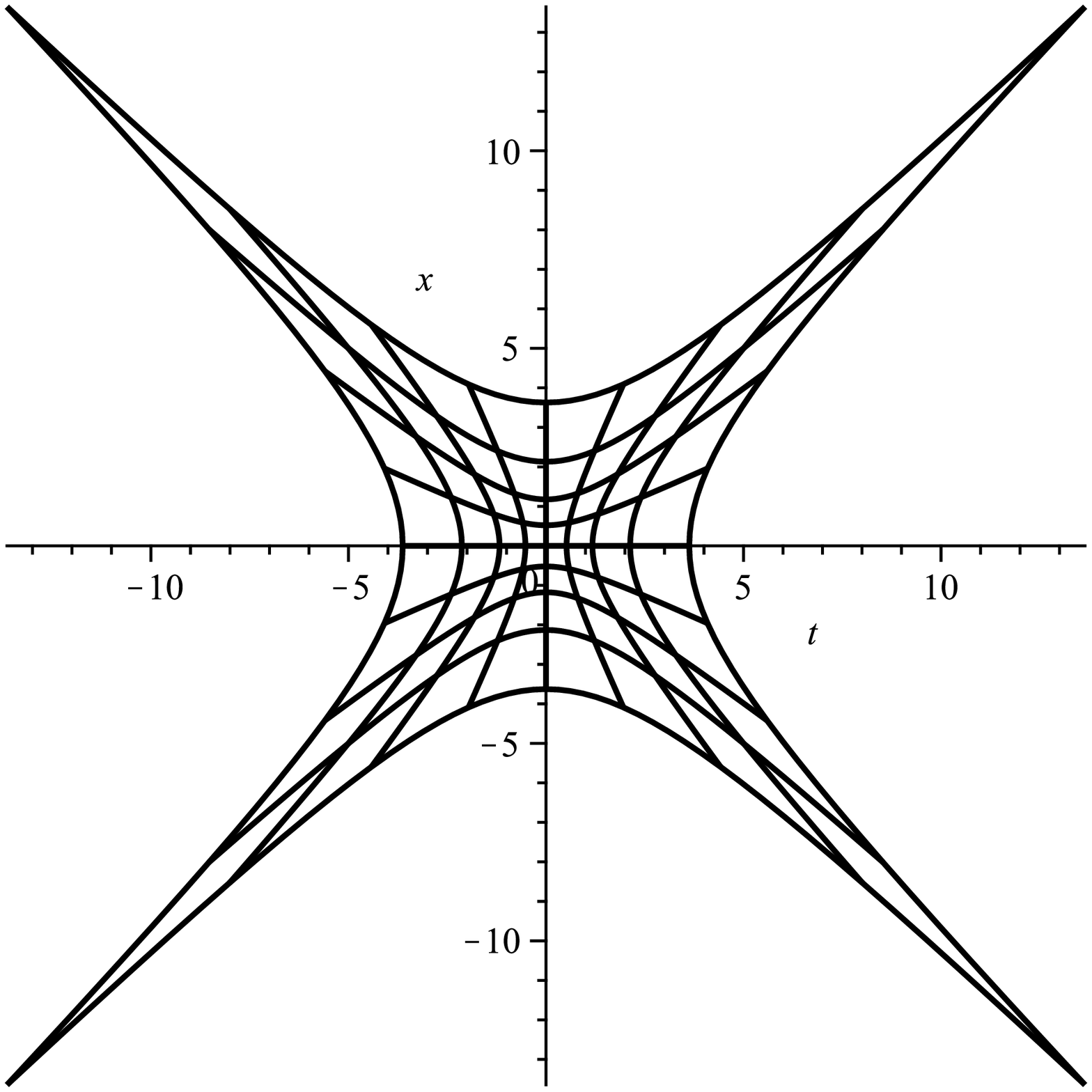}
    \medskip\nopagebreak\\ Fig.~\thefigure. The structure of the mapping $h\mapsto \sinh h$. }\par\medskip

The mapping $h\mapsto\sinh h$ is bijective, and hence its inverse mapping $\text{Arsh}$
    is defined on the whole double plane. Its explicit expression in coordinates is given by the formula:
\[\text{Arsh}\,h=\frac{1}{2}\left(\text{Arsh}[(t+x)\sqrt{1+(t-x)^2}+(t-x)\sqrt{1+(t+x)^2}]+\right.\]
    \[\left.j\text{Arsh}[(t+x)\sqrt{1+(t-x)^2}-(t-x)\sqrt{1+(t+x)^2}]\right).\]
%
%

Due to the duality of the hyperbolic cosine, the mapping
\[\cosh h=\cosh t\cosh x-j\sinh t\sinh x\]
is of a different kind. 
%

The first quadrant with its vertex in zero is bijectively mapped by $\cosh$ into
    the first quadrant with the vertex at the point $1$. Under these circumstances,
    the Cartesian net is mapped into a net of orthogonal hyperbolas.
The remaining quadrants with the vertex at the point $0$ are mapped as well to this quadrant (Fig.\ref{hypc}).

{\centering\small\refstepcounter{figure}\label{hypc}
    \includegraphics[width=.4\textwidth]{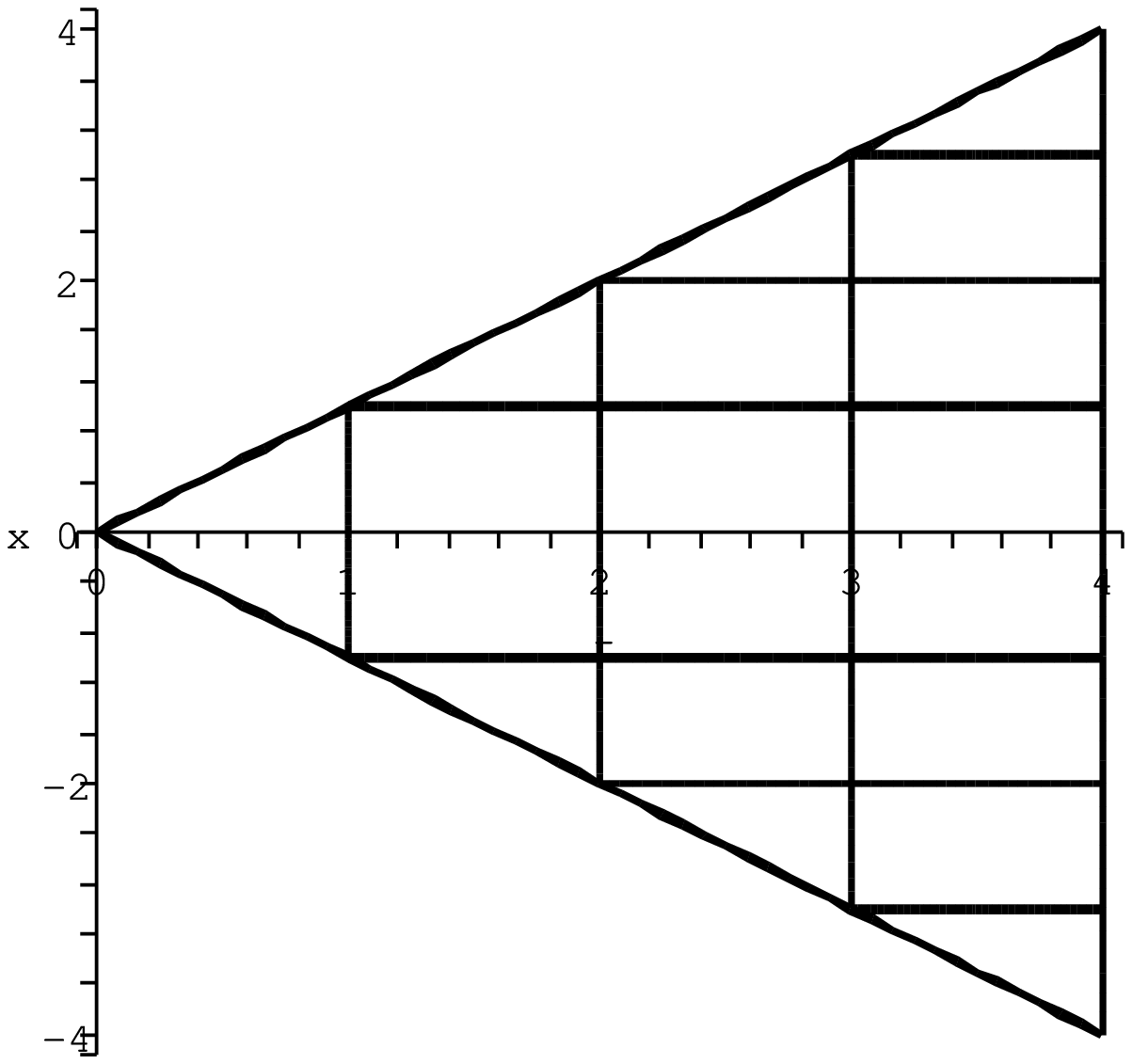}\includegraphics[width=.4\textwidth]{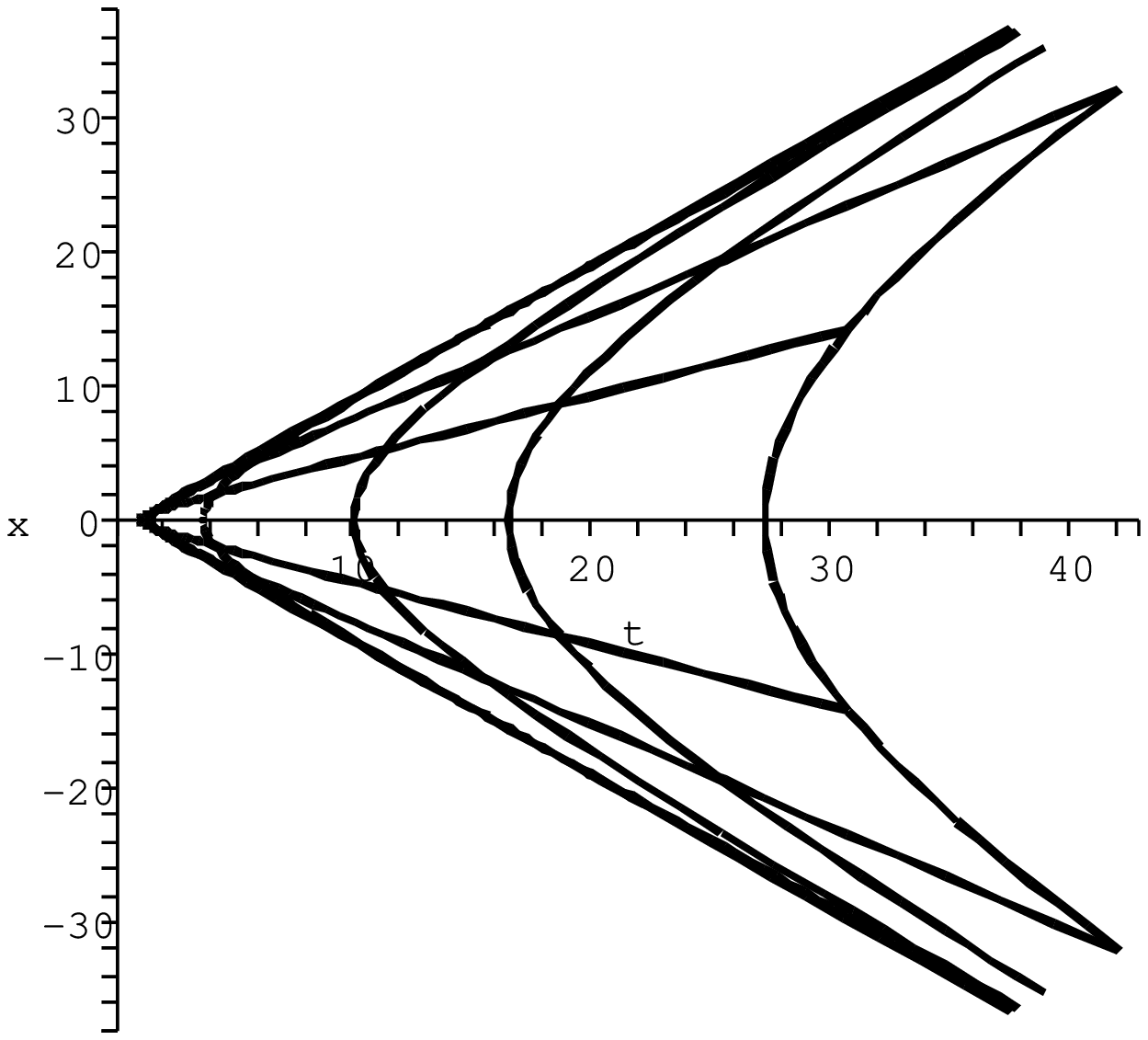}
    \medskip\nopagebreak\\ Fig.~\thefigure. The structure of the mapping $h\mapsto \cosh h$. }\par\medskip

The global structure of the mapping $\cosh h$ is illustrated in Fig.\ref{quad},
    from which the hashed quadrant is shifted to the right with one unit. In this way,
    the hyperbolic cosine is a 4-fold mapping, and the hyperbolic arccosine is 4-valent, with
    its Riemannian surface, represented in Fig.\ref{riemind}. Its explicit coordinate
    representation is given by the formula:
\[\text{Arch}\,h=\frac{1}{2}(\text{Arch}[t^2-x^2-\sqrt{(t+x)^2-1}\sqrt{(t-x)^2-1}]+\]
\[j\text{Arch}[t^2-x^2+\sqrt{(t+x)^2-1}\sqrt{(t-x)^2-1}]).\]
%
%

The function
\[\tanh h\equiv\frac{\sinh h}{\cosh h}=\frac{\tanh t(1-\tanh^2x)}{1-\tanh^2 t\tanh^2 x}+
    j\frac{\tanh x(1-\tanh^2 t)}{1-\tanh^2 t\tanh^2 x}\]
maps the double plane into the interior of the square with
vertexes  $(1,0),$ $(0,1),$ $(-1,0),$ $0,-1,$ and the square with
side 2 --- into some interior  domain near to the square-image of
double plane (Fig.\ref{tanh}).

{\centering\small\refstepcounter{figure}\label{tanh}
    \includegraphics[width=.4\textwidth]{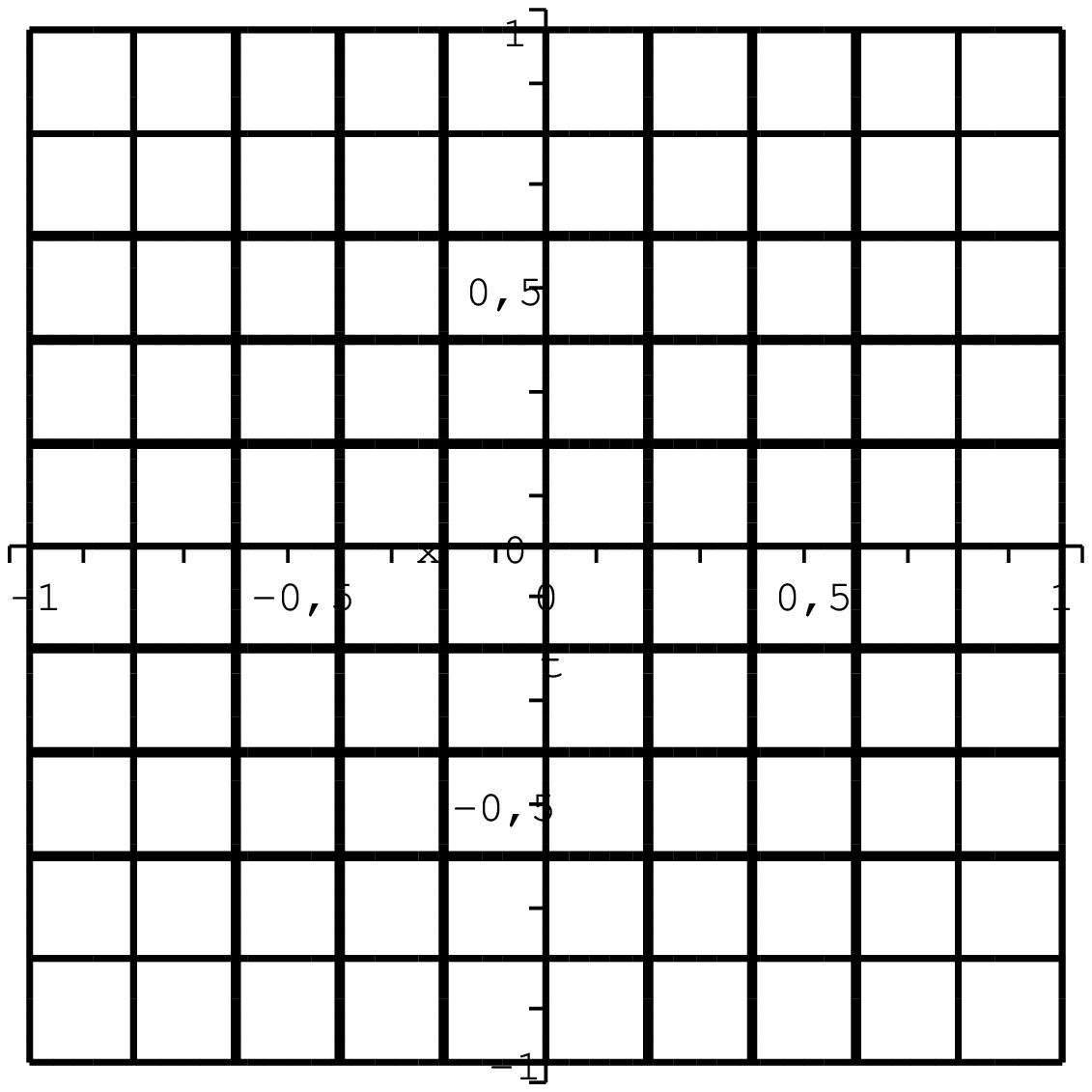}\includegraphics[width=.4\textwidth]{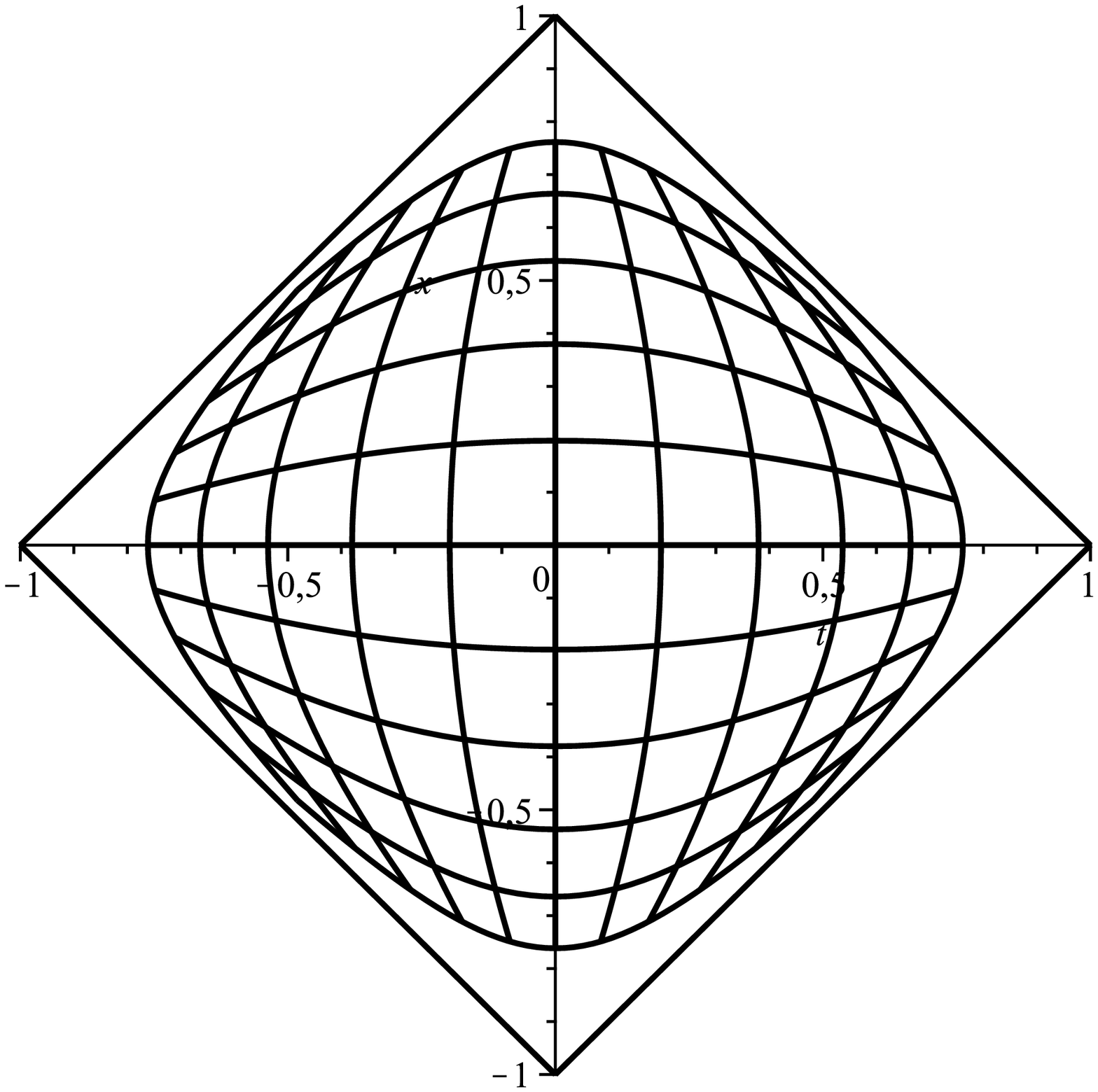}
    \medskip\nopagebreak\\ Fig.~\thefigure. The structure of the mapping $h\mapsto \tanh h$. }\par\medskip

%

The hyperbolic tangent is a mapping with one leaf, and hence its inverse mapping $\text{Arth}$
    is bijective on its domain. Its coordinate expression is given by the formula:
\[\text{Arth}\,h=\frac{1}{2}(\text{Arth}\,(t+x)+\text{Arth}\,(t-x)+j(\text{Arth}\,(t+x)-\text{Arth}\,(t-x))).\]
%
%

The function
\[\coth h\equiv \frac{\cosh h}{\sinh h}=\frac{\coth t(1-\coth^2x)+j\coth x(1+\coth^2 t)}{\coth^2t-\coth^2 x}\]
    is, in a certain sense, complementary to the function $\tanh h$: it maps the
    whole plane $\H $ into the exterior of the square in Fig.\ref{tanh} (at right).

Here, the rectilinear coordinate net is mapped to an orthogonal family of hyperbolas, which intersect
    at an infinity-point.
%

The function $\coth$ has one leaf and its inverse mapping $\text{Arcth}$ is univalent in its
    definition domain. Its coordinate expression is given by the formula:

\[\text{Arcth}\,h=\frac{1}{2}(\text{Arcth}\,(t+x)+\text{Arcth}\,(t-x)+j(\text{Arcth}\,
    (t+x)-\text{Arcth}\,(t-x))).\]

%

%
\subsection{The homographic mapping: $h\mapsto (ah+b)(ch+d)^{-1}$}

We define the homographic function by the relation:
\begin{equation}\label{drobno}w=\frac{ah+b}{ch+d}=D_{cd}^{ab}(h),\end{equation}
where $a,b,c,d$ are arbitrary double numbers which satisfy the condition $ad-bc\neq0$.
    The special character of the definition (\ref{drobno}) resides in the occurrence of
    the cone $\text{Con}(-d/c)$ on which this mapping is not well defined.
%

\begin{enumerate}\item
    The homographic transformation (\ref{drobno}) leads to a bijective and continuous
    $\H \setminus\text{Con}(-d/c)\to\H $ mapping in its domain of
    definition, which is conformal. Its inverse mapping is also homographic and has the form:
\begin{equation}\label{drobnoin}
    h=\frac{dw-b}{a-cw}=D^{d,-b}_{-c,a}(w)=(D_{cd}^{ab})^{-1}(h),\quad
    w\not\in\text{Con}(a/c).\end{equation}
\item
    The composition of two homographic mappings $D_2\circ D_1$ is also a homographic mapping.
    The set of all homographic mappings forms a group which is isomorphic to SL$(2,\H )$.
\item
    The homographic mapping transforms the hyperbolic circle $\text{HS}_R(h_0)$ given by the equation
    $|h-h_0|^2=\pm R^2$ into a hyperbolic circle.
\item
    By introducing the definition of points $h$ and $h^\vee$, which are conjugate relative to the
    circle $\text{HS}_R(h_0)$, similar to the corresponding definition on the complex plane,
    (conjugate points lie on the same ray which emerges from the center, one being located inside
    and the other outside the circle and their hyperbolic distances to the center satisfy the relation:
    $|h-h_0|\cdot|h^\vee-h_0|=R^2$), we can state another geometric property of homographic mappings:
    the points, which are conjugate relative to $\text{HS}_R(h_0)$, are mapped to points which are
    conjugate relative to the circle $D(\text{HS}_R(h_0))$.
\end{enumerate}


The proof of the latter assertions formally repeat the proof for the similar statements from the
    theory of homographic functions of complex variable.


{\centering\small\refstepcounter{figure}\label{drobb}
    \includegraphics[width=.4\textwidth]{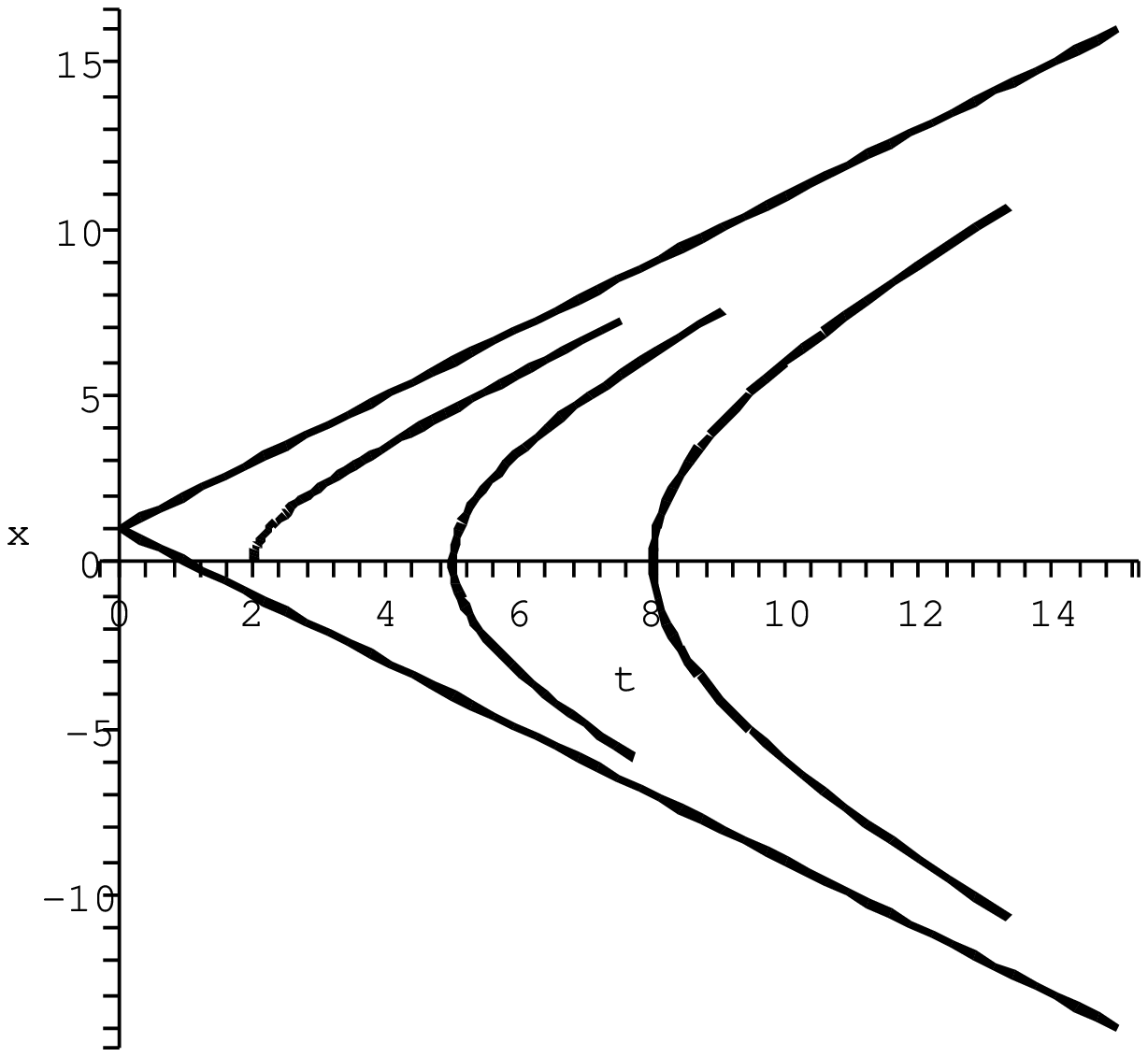}\includegraphics[width=.4\textwidth]{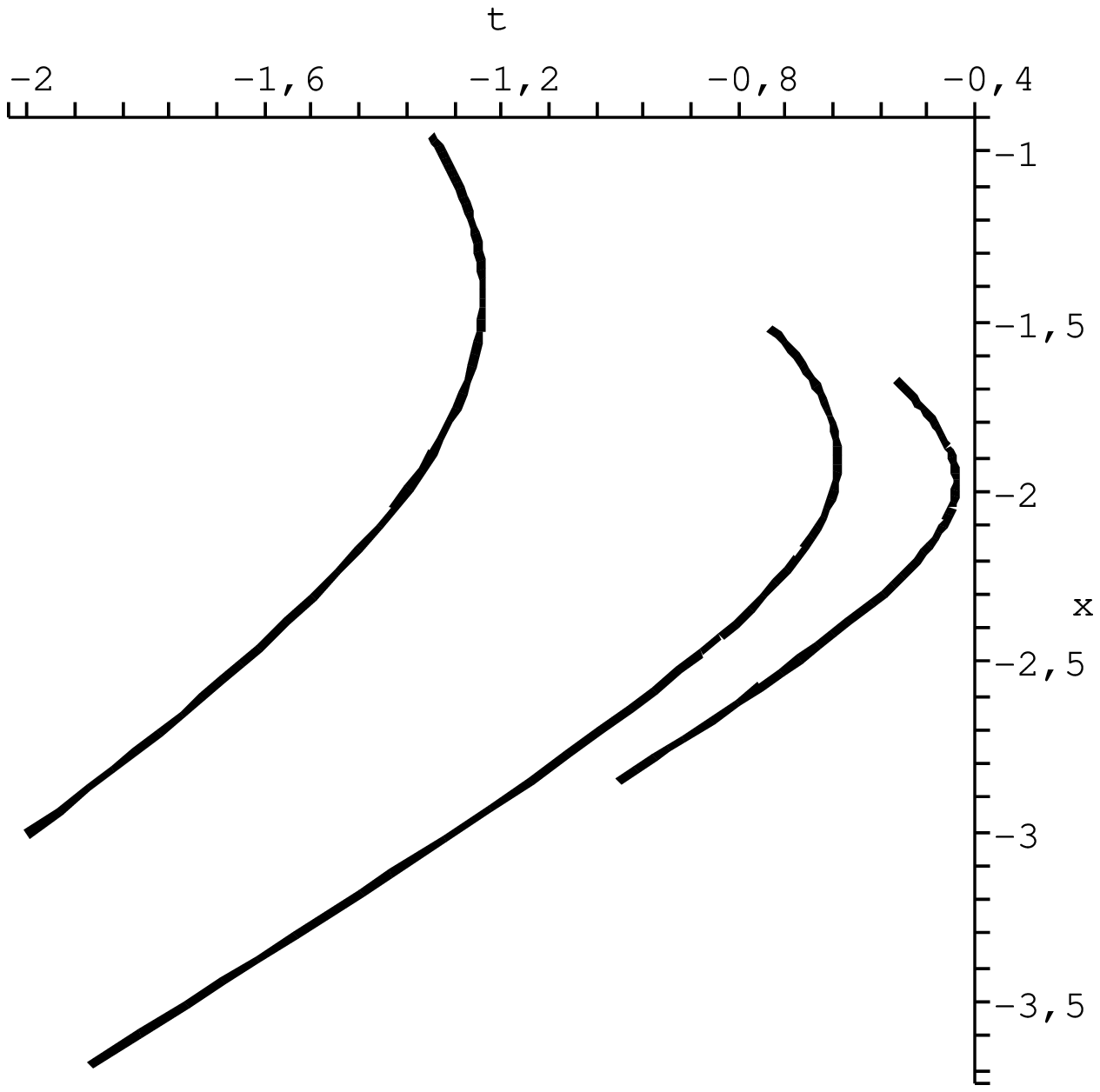}
    \medskip\nopagebreak\\ Fig.~\thefigure. The mapping of two arcs of hyperbolic circles into arcs of
    circles by means of the homographic mapping $h\mapsto (2h+j)/(1-j h)$. The pair of straight lines
    from the left picture represent the cone $\text{Con}(j)$}.\par\bigskip

%

\begin{figure}[htb]\centering\input{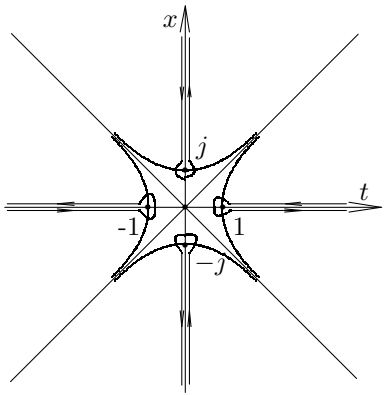}
    \caption{\small The action of the Zhukowskiy hyperbolic function. It maps components of  unit hyperbolic circle
    into corresponding parts of coordinate axes and twicely covers exteriority of square. }\label{zhuk}
\end{figure}

We shall examine, at last, the hyperbolic version of the Zhukowskiy function:
\[Z(h)\equiv\frac{1}{2}\left(h+\frac{1}{h}\right).\]
This transformation twicely maps double plane onto the plane with
removed unit square such as depicted in Fig. \ref{zhuk}.

Like in the case of complex variable, the Zhukowskiy mapping has two leaves:
    the outer part of the unit square and its inner part, which it maps onto the exteriority of the such square.
    At points $\pm1,\pm j$ the conformality of functions $Z(h)$ is broken, since at these points
    we have $Z'(h)=0$.




%
%
\section{Solving plane initial-boundary problems of the 2-dimensional field theory}

By analogy to the applications of conformal transformations on the standard (elliptic)
    complex plane, the $h$-analytic functions can be used for solving problems of field theory,
    which are related to the 2-order wave equation: $\Box_2\varphi=0$.
The lack of acknowledgement which the hyperbolic conformal mapping receive is firstly related
    by the non-traditional way in which are posed the initial-boundary problems, which can be solved
    by the method of hyperbolic conformal transformations.
Indeed, the use of conformal transformations in the plane of complex variable for solving
    problems of elliptic type, which are related to the Laplace equation, relies on the
    circumstances which were formerly described: the holomorphic mapping $w=f(z)$, whose real part
    represents the solution of some boundary problem, maps the boundary from outside the sources
    of the field into the straight line $\text{Re}\, w=\text{const}$.
This requirement reflects the condition of having a constant potential on the boundary of the domain,
    in which we solve the Laplace equation and guarantees the uniqueness of the solution up to
    an arbitrary choice of the value of potential on the boundary.

In initial-boundary problems of hyperbolic type we also use another setting of the problem:
    one examines the Space-Time domain as a multiple-boundary cylinder
    $D^3\times R_+$ (or a topologically equivalent to it alias)%
\footnote{In our case, $D^3$ is the 3-dimensional ball,
$R_+=[0,\infty)$},
    and there are provided the initial conditions on the initial surface $D^3\times\{0\}$
    (the initial values for the field and its derivatives by time) and the initial-boundary conditions
    on the flank surface $\partial D^3\times R_+$ (the values for the field and its derivatives
    by space coordinate).
If the problem is well-posed, then these initial boundary condition data lead a unique solution
    with good properties at each moment of time $t>0$. In the 2-dimensional Space-Time, the
    boundary of the domain represents a time-like rectangle $I\times R_+$ or a topologically equivalent
    to it figure.

The use of conformal transformations, which are $h$-holomorphic
mappings, assumes passing from
    plane elliptic problems into plane hyperbolic problems. In other words, an
    $h$-holomorphic
    mapping $w=f(h)$ represents a solution of some initial-boundary problem, namely the one,
    for which this mapping transforms the 1-dimensional border from the domain located outside sources
    into the line $\text{Re}\, w=\text{const}$.
It is obvious, that such a posing of the problem is different from the standard one, since
    the providing of the initial-boundary conditions change to {\em providing the form
    of the world-surface (lines) of constant potential}.
This surface has a Space-Time character. In principle, this can be obtained by means of measuring
    the wave field $\varphi$ at different points of the space at different moments of time
    by means of an appropriate, sufficiently large number of devices.
The event points of Space-Time, given by the readings of the devices provide $\varphi=\text{const}$
    and give the needed surface.
According to the above considerations, the form of this surface uniquely defines
    the solution of the wave equation. However, such a posing of the problem is seldom used in practice,
    since the data are distributed in Space-Time.

As an example, we examine the problem of determining the wave field, which has the constant value
    $\varphi_0$ on the hyperbolic circle $t^2-x^2=R^2$. According to the results of section \ref{exppp},
    the appropriate solution has the form: $\varphi=\varphi_0+\ln[(t^2-x^2)/R^2]=\text{Re}
    (\ln h+\varphi_0)$. Its 3-dimensional graph and the subsequent time slices are shown in Fig.\ref{3ddd}.

{\centering\small\refstepcounter{figure}\label{3ddd}
    \includegraphics[width=.35\textwidth]{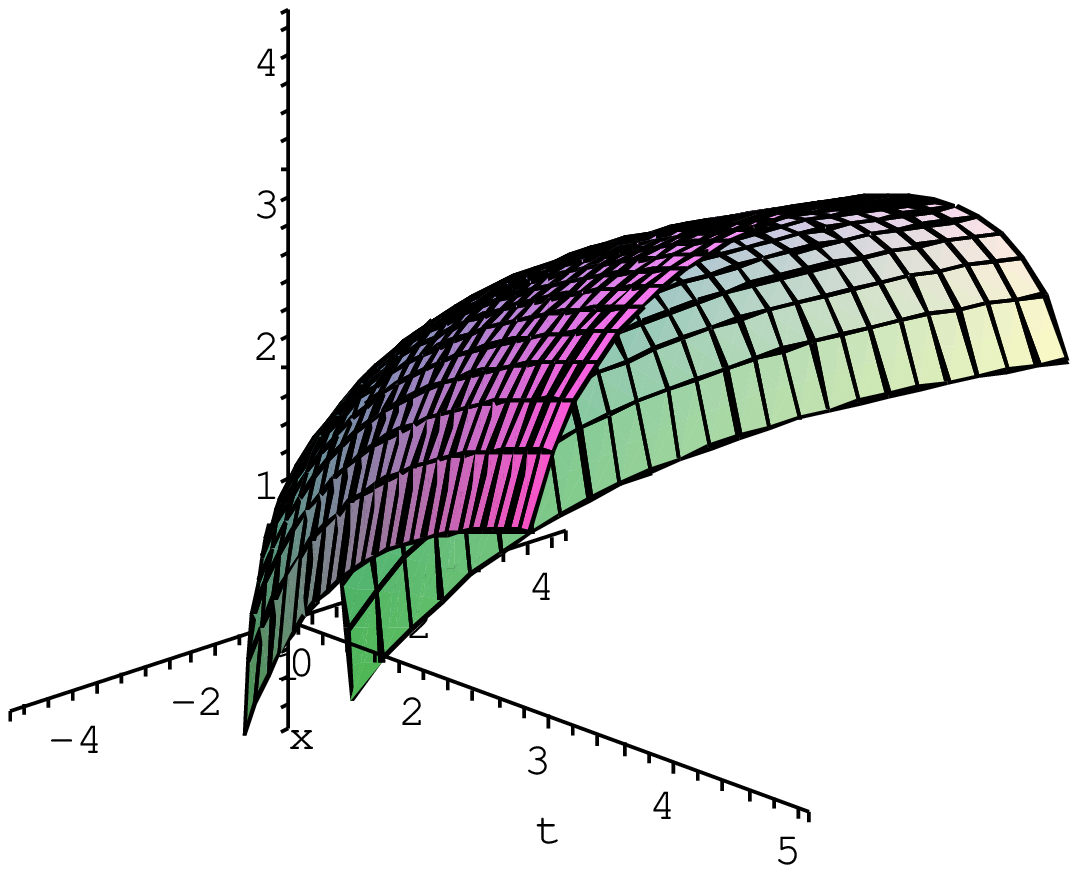}\includegraphics[width=.35\textwidth]{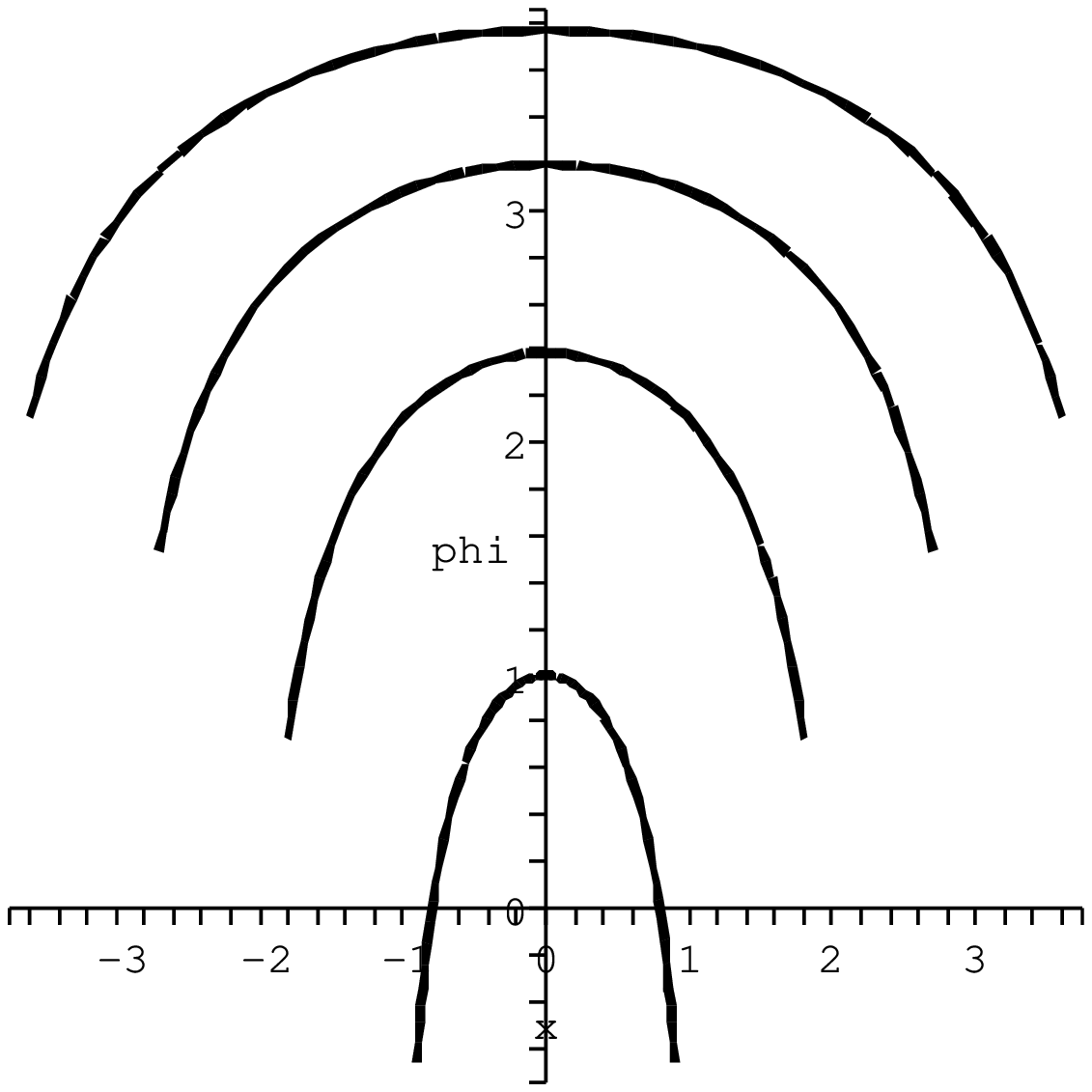}
    \medskip\nopagebreak\\ Fig.~\thefigure. The wave solution, which is constant on the
    hyperbolic circle $(R=1,\,\varphi_0=1)$. In the right picture there are represented
    the subsequent time slices of the surface, which is shown in the left picture, for $t=1,2,3,4$}.
%
%
\section{Conclusions}

We showed the general features of the fundamental for physics analogy between
    complex and double numbers. The algebraic description of the geometry and physics
    in 2-dimensional Space-Time by means of double numbers is not less appealing than
    describing the problems on the Euclidean plane by means of of complex numbers.
We have left non-addressed a series of problems (which essentially have a mathematical character),
    which appear along this path, giving priority to resemblance more than to rigor.
The difficulties encountered in the analysis based on the natural metric (i.e.,
    in the case examined by us, pseudo-Euclidean) topology are well known to the specialists in
    mathematical Relativity Theory and these (together with the difficulties of psychologic character)
    were the reason for rather low interest from physicists and mathematicians towards double
    numbers regarding their applications.
In our opinion, these difficulties can be in principle removed, and moreover,
    what is usually considered to be "difficulty"\, is in fact the expression
    of several new fundamental concepts, which the geometry of Space-Time brings to Physics and Mathematics.

The formerly described tight correspondence between the complex and the double planes
    can be also extended to vector fields. Like for each holomorphic function on the complex plane
    there corresponds some potential or solenoidal vector field, related to sources, vortices,
    dipoles, quadrupoles, etc, at the points where the holomorphicity of a $h$-holomorphic function is lost,
    we can also give sense to some vector field in the 2-dimensional Space-Time, whose singularities
    now have not elliptic, but hyperbolic properties.
Thus, we have the possibility to talk about hyperbolic positive and negative sources,
    hyperbolic vortices and vortex-wells, multipoles, Zhukowskiy hyperbolic function, etc.
The vector fields which are generated by these $h$-holomorphic
mappings have remarkable properties,
    like potentiality and solenoidality, which are obviously apprehended in a hyperbolic meaning.


Another interesting consequence provided by the $h$-holomorphic
functions of double variable is related to
    the possibility of regarding their associated conformal mappings, as changes from non-inertial
    reference systems in the flat 2-dimensional Space-Time.
Then the study of the group of conformal symmetries naturally extends the frames of 2-dimensional
    Special Relativity Theory, in which there are usually considered only the isometric mappings,
    regarded as changes between inertial reference systems.
Such a possibility (due to the absence of a corresponding infinite-dimensional conformal symmetry group)
    does not exist in three and four-dimensional pseudo-Euclidean spaces, but still, this exists in
    two dimensions and also in several three and four-dimensional flat Finsler geometries, e.g.,
    in spaces with Berwald-Moor metric, whose particular case is, in fact, the space of double numbers
    itself (\cite{4,3}).
The idea of a fundamental symmetry group (linear isometric, and non-linear conformal) was widespread
    in physics, though the use of the infinite conformal groups was in its essence limited to
    the case of 2-dimensional geometry with quadratic metrics.
While passing from quadratic metric forms of Space-Time to $n$-ary metric forms, this limitation
    is sometimes removed, and the problem of Finslerian extensions of SRT
    and GRT becomes actual and promising.
We need to emphasize that we talk about extensions, and not about any rejection of the
    classic representations, since we assume that the symmetry group which lie at the basis of
    previous physical-mathematical constructions, appear as in some sense secondary  (induced) fact
    of  the new infinite-dimensional symmetry groups existence.
Due to the limited space offered by a journal article, the authors plan to develop the present topic
    in a forthcoming paper.\par\medskip
%
%
{\bf Acknowledgement.} The authors express their gratitude to the researchers from NII-GKSGF V.~M. Chernov
and G. I. Garas'ko for examining the text and for their valuable remarks.
{\small\noindent
D.G. Pavlov\\
Research Institute of Hypercomplex Systems in Geometry and Physics,\\
Fryazino, Russia.\\
E-mail: geom2004@mail.ru\\\\
S.S. Kokarev\\
RSEC "Logos", Yaroslavl, Russia,\\
E-mail: logos-center@mail.ru
}
\end{document}